\documentclass[12pt,amscd]{amsart}
\usepackage{amssymb}

\numberwithin{equation}{section}
\newtheorem{thm}[equation]{Theorem}

\newtheorem{cor}[equation]{Corollary}

\newtheorem{lem}[equation]{Lemma}

\newtheorem{prop}[equation]{Proposition}

\newenvironment{pf}{\proof[\proofname]}{\endproof}
\newenvironment{pf*}[1]{\proof[#1]}{\endproof}

\theoremstyle{definition}

\newtheorem{defn}[equation]{Definition}

\theoremstyle{remark}

\newtheorem*{rmk}{Remark}
\newtheorem*{rmks}{Remarks}

\newtheorem*{ack}{Acknowledgement}

\newcommand{\comment}[1]{}

\makeatother

\begin{document}
\baselineskip=18truept


\def\C {{\mathbb C}}
\def\Cn {{\mathbb C}^n}
\def\R {{\mathbb R}}
\def\Rn {{\mathbb R}^n}
\def\Z {{\mathbb Z}}
\def\N {{\mathbb N}}
\def\cal#1{{\mathcal #1}}
\def\bb#1{{\mathbb #1}}

\def\dbar {\bar \partial }
\def\dir {{\mathcal D}}
\def\lev#1{{\mathcal L}(#1)}
\def\lap {\Delta }
\def\ol {{\mathcal O}}
\def\E {{\mathcal E}}
\def\J {{\mathcal J}}
\def\U {{\mathcal U}}
\def\V {{\mathcal V}}
\def\z {\zeta }
\def\Harm {\text {Harm}\, }
\def\grad {\nabla }
\def\dexh {\{ M_k \} _{k=0}^{\infty } }
\def\sing#1{#1_{\text {sing}}}
\def\reg#1{#1_{\text {reg}}}

\def\lquotient{\big{\backslash}}

\def\setof#1#2{\{ \, #1 \mid #2 \, \} }

\def\image#1{\text{\rm im}\,\bigl[#1\bigr]}
\def\kernel#1{\text{\rm ker}\,\bigl[#1\bigr]}

\def\holecl {M\setminus \overline M_0}
\def\hole {M\setminus M_0}

\def\nd{\frac {\partial }{\partial\nu } }
\def\ndof#1{\frac {\partial#1}{\partial\nu } }

\def\pdof#1#2{\frac {\partial#1}{\partial#2}}

\def\cinf{\ensuremath{C^{\infty}} }
\def\cinfns{\ensuremath{C^{\infty}}}

\def\diam{\text {\rm diam} \, }

\def\real{\text {\rm Re}\, }

\def\imag{\text {\rm Im}\, }

\def\supp{\text {\rm supp}\, }

\def\Vol{\text {\rm vol} \, }

\def\restrict#1{\upharpoonright_{#1}}

\def\sm{\setminus}

\def\plshclass{{\mathcal {P}}}
\def\strplshclass{{\mathcal S\mathcal P}}

\def\dist{{\text{\rm dist}}\,}
\def\distg#1{{\text{\rm dist}_{#1}}}

\def\geqtrace#1#2{{\geq_{(#1,#2)}}}
\def\gtrace#1#2{{>_{(#1,#2)}}}
\def\leqtrace#1#2{{\leq_{(#1,#2)}}}
\def\ltrace#1#2{{<_{(#1,#2)}}}



\def\anal{analytic }
\def\analns{analytic}

\def\bdd{bounded }
\def\bddns{bounded}

\def\cpt{compact }
\def\cptns{compact}

\def\cpx{complex }
\def\cpxns{complex}

\def\cont{continuous }
\def\contns{continuous}

\def\dime{dimension }
\def\dimens{dimension }

\def\exh{exhaustion }
\def\exhns{exhaustion}

\def\fn{function }
\def\fnns{function}

\def\fns{functions }
\def\fnsns{functions}

\def\holo{holomorphic }
\def\holons{holomorphic}

\def\mero{meromorphic }
\def\merons{meromorphic}

\def\holoconvex{holomorphically convex }
\def\holoconvexns{holomorphically convex}

\def\ircomp{irreducible component }
\def\concomp{connected component }
\def\ircompns{irreducible component}
\def\concompns{connected component}
\def\ircomps{irreducible components }
\def\concomps{connected components }
\def\ircompsns{irreducible components}
\def\concompsns{connected components}

\def\irred{irreducible }
\def\irredns{irreducible}

\def\con{connected }
\def\conns{connected}

\def\comp{component }
\def\compns{component}
\def\comps{components }
\def\compsns{components}

\def\mfld{manifold }
\def\mfldns{manifold}
\def\mflds{manifolds }
\def\mfldsns{manifolds}

\def\nbd{neighborhood }
\def\nbds{neighborhoods }
\def\nbdns{neighborhood}
\def\nbdsns{neighborhoods}

\def\harm{harmonic }
\def\harmns{harmonic}
\def\plh{pluriharmonic }
\def\plhns{pluriharmonic}
\def\plsh{plurisubharmonic }
\def\plshns{plurisubharmonic}

\def\qplsh#1{$#1$-plurisubharmonic}
\def\hplsh{$(n-1)$-plurisubharmonic }
\def\hplshns{$(n-1)$-plurisubharmonic}

\def\para{parabolic }
\def\parans{parabolic}

\def\rel{relatively }
\def\relns{relatively}

\def\str{strictly }
\def\strns{strictly}

\def\strg{strongly }
\def\strgns{strongly}

\def\cvx{convex }
\def\cvxns{convex}

\def\wrt{with respect to }
\def\wrtns{with respect to}

\def\st {such that }
\def\stns {such that}

\def\hm {harmonic measure }
\def\hmns {harmonic measure}

\def\hmib {harmonic measure of the ideal boundary of }
\def\hmibns {harmonic measure of the ideal boundary of}

\def\Vert{\text{{\rm Vert}}\, }
\def\Edge{\text{{\rm Edge}}\, }

\def\til#1{\tilde{#1}}
\def\wtil#1{\widetilde{#1}}

\def\what#1{\widehat{#1}}

\def\seq#1#2{\{#1_{#2}\} }


\def\vphi {\varphi }


\def\inv{   ^{-1}  }

\def\ssp#1{^{(#1)}}

\def\set#1{\{ #1 \}}

\title[$L^2$ Castelnuovo-de Franchis]
{$L^2$ Castelnuovo-de Franchis, the cup product lemma, and
filtered ends of K\"ahler manifolds}
\author[T.~Napier]{Terrence Napier$^{*}$}
\address{Department of Mathematics\\Lehigh University\\Bethlehem, PA 18015}
\email{tjn2@lehigh.edu}
\thanks{$^{*}$Research partially
supported by NSF grant DMS0306441}
\author[M.~Ramachandran]{Mohan Ramachandran}
\address{Department of Mathematics\\SUNY at Buffalo\\Buffalo, NY 14260}
\email{ramac-m@math.buffalo.edu}


\date{November 6, 2007}

\begin{abstract}
Simple approaches to the proofs of the $L^2$ Castelnuovo-de~Franchis
theorem and the cup product lemma which give new versions are developed.
For example, suppose $\omega_1$
and $\omega_2$ are two linearly independent closed \holo $1$-forms on a bounded geometry \con
complete K\"ahler manifold~$X$ with $\omega_2$ in~$L^2$. According to a version of the
$L^2$ Castelnuovo-de~Franchis theorem obtained in this paper,
if $\omega_1\wedge\omega_2\equiv 0$, then there exists a surjective
proper \holo mapping of~$X$ onto a Riemann surface for which $\omega_1$ and $\omega_2$
are pull-backs.  Previous versions required both forms to be in~$L^2$.
\end{abstract}

\maketitle

\section*{Introduction} \label{introduction}

According to the classical theorem of Castelnuovo and de~Franchis
(see \cite{Be}, \cite{BPV}), if, on a \con \cpt \cpx manifold~$X$,
there exist linearly independent closed \holo $1$-forms $\omega_1$
and $\omega_2$  with $\omega _1\wedge \omega _2\equiv 0$, then
there exist a surjective \holo mapping $\Phi$ of $X$ onto a
curve~$C$ of genus~$g\geq 2$ and \holo $1$-forms $\theta_1$ and
$\theta_2$ on~$C$ \st $\omega_j=\Phi^*\theta_j$ for $j=1,2$. The
main point is that the meromorphic \fn $f\equiv\omega_1/\omega_2$
actually has no points of indeterminacy, so one may Stein factor
the \holo map $f\colon X\to\mathbb P^1$.
\begin{rmk}
The requirement that the forms be closed is superfluous if the
\cpt manifold $X$ is a surface or if $X$ is K\"ahler. For, if
$\eta=\sum \sqrt{-1}g_{i\bar j}dz_i\wedge d\bar z_j$ is the
K\"ahler form for a K\"ahler metric~$g$ and $\omega$ is a
\holo $1$-form, then, by Stokes' theorem, we have
\[
\int_Xd\omega\wedge d\bar\omega\wedge\eta^{n-2}=0;
\]
where $n=\dim X$. Since the integrand is a nonnegative $2n$-form,
the form must vanish and it follows that $d\omega=0$. For $X$ a
surface, the same argument with the factor $\eta^{n-2}$ removed
again yields $d\omega=0$.
\end{rmk}

In general, given a \con \cpx manifold $X$ and linearly
independent closed \holo $1$-forms $\omega_1$ and $\omega_2$ on
$X$ with $\omega_1\wedge\omega_2\equiv 0$, the \mero \fn
$f\equiv\omega_1/\omega_2$ has no points of indeterminacy, $f$ is
locally constant on the \anal set $Z=\setof{x\in
X}{(\omega_1)_x=0\text{ or }(\omega_2)_x=0}$, and $f$ is constant
on each leaf of the \holo foliation determined by $\omega_1$ and
$\omega_2$ in $X\sm Z$ (see, for example, \cite{NR-BH Weakly
1-complete} for an elementary proof). In particular, if the levels
of the \holo map $f\colon X\to\mathbb P^1$ are \cptns, then Stein
factorization gives a surjective proper \holo mapping of $X$ onto
a Riemann surface.

We will say that a complete Hermitian manifold $(X,g)$ has
\textit{bounded geometry of order~$k$} if, for some constant $C>0$
and for every point $p\in X$, there is a biholomorphism $\Psi $ of
the unit ball $B=B(0;1)\subset \C ^n$ onto a \nbd of $p$ in $X$
\st $\Psi (0)=p$ and, on $B$,
$$ C\inv
g_{\C ^n} \leq \Psi ^*g\leq Cg_{\C ^n} \quad\text{and}\quad |
D^m\Psi ^*g| \leq C \text{ for } m=0,1,2,\dots,k.
$$
For $k=0$, we will simply say that $(X,g)$ has \textit{bounded
geometry}.
Gromov~\cite{Gro-Kahler hyperbolicity} observed
that, for $f=\omega_1/\omega_2$ as above, one gets \cpt levels if
$X$ is a bounded geometry complete K\"ahler manifold and the
$1$-forms are in~$L^2$ and have exact real parts; thus giving an
$L^2$ version of the Castelnuovo-de~Franchis theorem. He also
introduced his so-called \textit{cup product lemma}, according to
which, two $L^2$ \holo $1$-forms $\omega_1$ and $\omega_2$ with
exact real parts on a bounded geometry complete K\"ahler manifold
must satisfy $\omega_1\wedge\omega_2\equiv 0$. He applied these
results to the study of K\"ahler groups. Other versions have since
been developed and applied by others in many different contexts.
Other versions and applications of the Castelenuovo-de~Franchis
theorem (for \cpt and non\cpt manifolds) and the cup product lemma
appear in, for example, \cite{Siu-strong rigidity}, \cite{CarT},
\cite{Gro-Sur la groupe fond}, \cite{Li Structure complete
Kahler}, \cite{Gro-Kahler hyperbolicity}, the work of Beauville
(see \cite{Cat1}), \cite{Simpson-VHS}, \cite{Gromov-Schoen},
\cite{ArapuraBressRam}, \cite{JY1}, \cite{JY2},
\cite{Simpson-Lefschetz thm integral leaves}, \cite{Arapura},
\cite{NR-Structure theorems}, \cite{ABCKT}, \cite{Mok-Castelnuovo
de Franchis Unitary}, \cite{Jost-Zuo}, \cite{NR-BH Weakly
1-complete}, \cite{NR-BH Regular hyperbolic Kahler},
\cite{Delzant-Gromov Cuts}, \cite{NR Thompson not Kahler}, and
\cite{NR Filtered ends}. In this paper, new approaches to the
proofs of the $L^2$~Castelnuovo-de~Franchis theorem and to the cup
product lemma are developed. These new approaches are simpler than
previous approaches and give more general results. In particular, a
version of the $L^2$~Castelnuovo-de~Franchis theorem is obtained in
which only \emph{one} of the \holo $1$-forms need be
in~$L^2$.

\begin{thm}[$L^2$ Castelnuovo-de~Franchis theorem]\label{L2 Castelnuovo de Franchis thm from intro}
Let $(X,g)$ be a \con complete K\"ahler manifold with bounded
geometry and let $\omega_1$ and $\omega_2$ be linearly independent
closed \holo $1$-forms on $X$ \st $\omega_1$ is in~$L^2$ and
$\omega_1\wedge\omega_2\equiv 0$. Then there exist a surjective
proper \holo mapping $\Phi\colon X\to S$ of $X$ onto a Riemann
surface~$S$ with $\Phi_*\ol_X=\ol_S$ and \holo $1$-forms
$\theta_1$ and $\theta_2$ on $S$ \st $\omega_j=\Phi^*\theta_j$ for
$j=1,2$.
\end{thm}

The main point of the proof is that, for a suitable small open set, the
holonomy induced by the \holo foliation associated to the \holo
$1$-forms is trivial (see Section~\ref{L2 Castel de
Franchis sect}). A version for a bounded geometry (of order~$2$) end is also obtained
(Theorem~\ref{L2 Castelnuovo for end thm}).

For the cup product lemma, the main point is that one obtains different
versions by considering positive forms rather than just \holo $1$-forms;
an observation which has its roots in the theory of currents and
which has been applied in other contexts to obtain related results.
Simple Stokes theorem
arguments together with Gromov's arguments then give myriad
versions of which only a few will be considered in this paper (see
Sections~\ref{cup product lemma section}~and~\ref{Cup product
lemma general sect}). For example, there is the following version in which
one of the forms is assumed to be in~$L^\infty$ instead of in~$L^2$ and the other form
need not have exact real part:
\begin{thm}\label{Cup product lemma general thm intro}
Let $\omega_1$ and $\omega_2$ be closed \holo $1$-forms on a \con
complete K\"ahler manifold~$X$ \st $\omega_1$ is bounded,
$\real(\omega_1)$ is exact, and $\omega_2$ is in~$L^2$. Then
\(\omega_1\wedge\omega_2=0\).
\end{thm}
\begin{rmk} By the Gaffney theorem~\cite{Gaffney}, an $L^2$ \holo $1$-form on
a complete K\"ahler manifold is automatically closed, so the
requirement that $\omega_2$ be closed is superfluous.
\end{rmk}

Theorem~\ref{L2 Castelnuovo de Franchis thm from intro}
and Theorem~\ref{Cup product lemma general thm intro} together
give the following:
\begin{cor}\label{Map to RS cor intro}
Let $\omega_1$ and $\omega_2$ be linearly independent closed \holo
$1$-forms on a \con complete K\"ahler manifold~$X$ with bounded
geometry \st $\omega_1$ is bounded, $\real(\omega_1)$ is exact,
and $\omega_2$ is in~$L^2$. Then there exist a surjective proper
\holo mapping $\Phi\colon X\to S$ of $X$ onto a Riemann
surface~$S$ and \holo $1$-forms $\theta_1$ and $\theta_2$ on $S$
\st $\omega_j=\Phi^*\theta_j$ for $j=1,2$.
\end{cor}
\begin{rmk}
Since an $L^2$ \holo $1$-form on a bounded geometry complete
K\"ahler manifold is bounded, the condition that $\omega_1$ is
bounded may be replaced with the condition that
$\omega_1$ is in $L^2$.
\end{rmk}

The proof of Theorem~\ref{L2 Castelnuovo de Franchis thm from
intro} appears in Section~\ref{L2 Castel de Franchis sect} and that of
Theorem~\ref{Cup product lemma general thm intro}
in Section~\ref{cup product lemma section}. As an application, the
results are shown in Sections~\ref{Filtered ends prelim
sect}~and~\ref{pf of filtered ends sect} to give a slightly
simplified proof of the main result of~\cite{NR Filtered ends}.
Further generalizations of the cup product lemma appear in
Section~\ref{Cup product lemma general sect}. Finally, a version
of the $L^2$ Castelnuovo-de~Franchis theorem for an end (which is
applied in \cite{NR-BH bound geom hyperbolic Kahler}) is proved in
Section~\ref{L2 Castelnuovo for an end sect}.

\begin{ack}  We would like to thank Domingo Toledo for useful
conversations.
\end{ack}

\section{Proof of the $L^2$ Castelnuovo-de Franchis
theorem}\label{L2 Castel de Franchis sect}

Given two linearly independent closed \holo $1$-forms $\omega_1$
and $\omega_2$ on a \con \cpx manifold $X$ with
$\omega_1\wedge\omega_2\equiv 0$, we get a nonconstant \holo map
\[
f=\frac{\omega_1}{\omega_2}\colon X\to\mathbb P^1.
\]
We have $(f_*)\wedge\omega_1=(f_*)\wedge\omega_2\equiv 0$ since,
on $f\inv(\C)=f\inv(\mathbb P^1\sm\set\infty)$,
\(df\wedge\omega_2=d\omega_1=0\). It follows that $f$ is locally
constant on the \anal set
\[
Z=\setof{x\in X}{(\omega_1)_x=0\text{ or }(\omega_2)_x=0}
\]
(in particular, $f(Z)$ is countable) and $f$ is constant on each
leaf of the \holo foliation determined by $\omega_1$ and
$\omega_2$ in $X\sm Z$. Thus $\omega_1$ and $\omega_2$ determine a
singular \holo foliation in~$X$ with closed leaves given by the
levels of~$f$. Moreover, for $j=1,2$, $\omega_j$ is exact in a
\nbd of each level~$L$ of $f$. For the integral of~$\omega_j$
along any closed loop in $L$ and, therefore, along any closed loop
in a small \nbd of $L$, must be zero.

The main step in the proof of Theorem~\ref{L2 Castelnuovo de
Franchis thm from intro} is the following:
\begin{lem}\label{Finite volume for level lemma} Let $(X,g)$ be a
\con Hermitian manifold. If $\omega_1$ and $\omega_2$ are two
linearly independent closed \holo $1$-forms on~$X$, $\omega_1$ is
in~$L^2$, $\omega_1\wedge\omega_2\equiv 0$, and \( f={\omega_1}/{\omega_2}\colon X\to\mathbb P^1\),
then the levels of~$f$ over almost every regular value have finite
volume (that is, almost every smooth (closed) leaf of the \holo
foliation determined by $\omega_1$ and $\omega_2$ has finite
volume). 
\end{lem}
\begin{pf}
Given a regular value $\zeta_0\in\C=\mathbb P^1\sm\set{\infty}$
of~$f$ and a point $p\in f\inv(\zeta_0)$, we may choose a \rel
\cpt \holo coordinate \nbd $(U,z=(z_1,\dots,z_n))$ in $X$ in which
$f\restrict U=z_1$, $p=(\zeta_0,0,\dots,0)$, and $U=D\times
\Delta^{n-1}$ where $D$ is a disk centered at $\zeta_0$ and $\Delta$ is
a disk centered at $0$ in $\C$; and we may choose a \holo \fn $h$
on $U$ with $\omega_1\restrict U=dh$.

If $A=D\times\set{0}\subset U$ and $\Omega$ is the union of
all of those levels of $f$ which meet~$A$, then $\Omega$ is a
nonempty \con open subset of $X$ containing~$U$. For if $\seq
x\nu$ is a sequence in $X$ converging to a point $y\in\Omega$, $L$
is the level containing~$y$, and $L_\nu$ is the level containing
$x_\nu$ for each $\nu$, then, by continuity of intersections (see
\cite{Stein}, \cite{Tworzewski-Winiarski Cont of intersect}, and
Section~4.3 of \cite{ABCKT}), after replacing the sequence with a
suitable subsequence, we get $L_\nu\to L$. Since $L$ meets $U$, we
have $L_\nu\cap U\neq\emptyset$, and hence $L_\nu\cap
A\neq\emptyset$, for $\nu\gg 0$ ($L_\nu\cap U$ and $L\cap U$ are
slices of the form $\set{\zeta}\times\Delta^{n-1}$). Thus
$x_\nu\in\Omega$ for $\nu\gg 0$ and it follows that $\Omega$ is
open.

Since $dh\wedge dz_1=\omega_1\wedge df\equiv 0$, $h$ is constant in the
variables $(z_2,\dots,z_n)$ in $U$ and we have $h=k(z_1)$ on $U$
for some nonconstant \holo \fn $k$ on $D$. Thus, since
$f(\Omega)=f(U)=D$, we may form the \holo extension
$h_0=k(f)$ on $\Omega$ and, since $dh=\omega_1$ on $U$, we get
$dh_0=\omega_1$ on $\Omega$. In particular, since $\omega_1$ is in
$L^2$, $h_0$ must have finite energy. Setting $u=\real(h_0)$ and
$v=\imag(h_0)$ and applying the coarea formula to the mapping
$(u,v)\colon\Omega\to\R^2$, we see that there exists a set~$S'$ of measure~$0$ in~$\C$
\st $\Vol(h_0\inv(\zeta))<\infty$ for each $\zeta\in \C\sm S'$. We may
choose a set~$S$ of measure~$0$ in~$D$ which contains the set of critical values
of~$f$ as well as the set $k\inv(S')$.
For each point $\zeta\in D\sm S$, the
level~$L=f\inv(\zeta)\cap\Omega$ of
$f$ over~$\zeta$ meeting~$\Omega$ is a \concomp of $h_0\inv(k(\zeta))$ and,
therefore, $\Vol (L)<\infty$.

Finally, forming a countable collection $\seq U\nu$ of such open
sets $U$ in $X$ covering
\[
f\inv\left(\set{\text{\,regular values\,}}\sm\infty\right),
\]
forming the associated measure~$0$ sets $\seq S\nu$
in~$\C\subset\mathbb P^1$, and letting $S\subset\mathbb P^1$ be
the measure~$0$ set given by
\[
S=\bigcup_\nu S_\nu\cup\set{\text{\,critical values\,}}
\cup\set\infty,
\]
we see that each of the levels of~$f$ over every point in $\mathbb P^1\sm
S$ has finite volume.
\end{pf}

Theorem~\ref{L2 Castelnuovo de Franchis thm from intro} now
follows from standard arguments (see \cite{Gro-Kahler
hyperbolicity}, \cite{ArapuraBressRam}, and Chapter~4
of~\cite{ABCKT}) which are sketched below for the convenience of
the reader.
\begin{pf*}{Proof of Theorem~\ref{L2 Castelnuovo de Franchis thm from
intro}} Let $(X,g)$ be a \con complete K\"ahler manifold with
bounded geometry and let $\omega_1$~and~$\omega_2$ be two linearly
independent closed \holo $1$-forms \st $\omega_1$ is in~$L^2$ and
$\omega_1\wedge\omega_2\equiv 0$. We may also assume that $n=\dim
X>1$. The \holo map $f=\omega_1/\omega_2\colon X\to\mathbb P^1$ is
open, and, by Lemma~\ref{Finite volume for level lemma}, we may
fix a regular value $\zeta_0\in f(X)\sm\set\infty$ and a \concomp
$L_0$ of the submanifold $f\inv(\zeta_0)$ of $X$ \st $\Vol
(L_0)<\infty $. Lelong's monotonicity formula (see 15.3 in
\cite{Chirka}) shows that there is a constant $c>0$ \st each point
$p\in X$ has a \nbd $U_p$ \st $\diam(U_p)<1$ and $\Vol(A\cap
U_p)\geq c$ for every \cpx \anal set $A$ of pure dimension $n-1$
in $X$ with $p\in A$. Therefore, since $L_0$ has finite volume,
$L_0$ must be \cptns.

It follows that the set $V=\setof{x\in X}{x\text{ lies in a \cpt
level of }f}$ is a nonempty open set. To show that $V$ is also
closed, let $V_0$ be a \comp of $V$, let $\seq xj$ be a sequence
in $V_0$ converging to a point $p\in \overline V_0$, and, for each
$j$, let $L_j\subset V_0$ be the \cpt level of $f$ through $x_j$.
Stein factoring $f\restrict{V_0}$, we get a proper \holo mapping
$\Phi\colon V_0 \to S$ onto a Riemann surface~$S$ with
$\Phi_*\ol_{V_0}=\ol_S$. We may choose each $x_j$ to lie over a
regular value of $f$ and of $\Phi$. Applying Stokes' theorem as in
\cite{Stoll}, we see that $\Vol (L_j)$ is constant in $j$ and so
the above volume estimate implies that, for some $R\gg 0$, we have
$L_j\subset B(p;R)$ for $j=1,2,3,\dots$. On the other hand, by
\cite{Stein} (see also \cite{Tworzewski-Winiarski Cont of intersect}
or Theorem~4.23 in \cite{ABCKT}), a subsequence of $\seq Lj$
converges to the level $L$ of $f$ through $p$. So we must have
$L\subset\overline {B(p;R)}$ and hence $L$ is \cptns. Thus $p\in\overline
V_0\cap V$ and, therefore, $p\in V_0$.
It follows that $V=V_0=X$. Thus every level of $f$ is \cpt and we
get our proper \holo mapping $\Phi\colon X\to S$.

Finally, we recall that, for each $j=1,2$, $\omega_j$ is exact on
a \nbd of each level of~$f$; that is,~on a \nbd of each fiber
of~$\Phi$. Thus, for each point $s\in S$, we have a \con \nbd $D$
of $s$ in $S$ and a \holo \fn $h_j$ on $U=\Phi\inv(D)$ \st
$\omega_j=dh_j$ on $U$. The \fn $h_j$ descends to a unique \holo
\fn $k_j$ on $D$ with $\Phi^*k_j=h_j$. Thus we get a unique
well-defined \holo $1$-form $\theta_j$ on $S$ with
$\Phi^*\theta_j=\omega_j$ by setting $\theta_j\restrict D=dk_j$ on
each such \nbd $D$.
\end{pf*}

The following easy consequence is a more convenient form for some
applications:
\begin{cor}\label{L2 Castelnuovo for partial exact forms cor}
Let $(X,g)$ be a \con complete K\"ahler manifold with bounded
geometry and let $\rho_1$ and $\rho_2$ be two real-valued \plh
\fns on $X$ \st $d\rho_1$ and $d\rho_2$ are linearly independent,
$\rho_1$ has finite energy, and
$\partial\rho_1\wedge\partial\rho_2\equiv 0$. Then there exist a
surjective proper \holo mapping $\Phi\colon X\to S$ of $X$ onto a
Riemann surface~$S$ with $\Phi_*\ol_X=\ol_S$ and real-valued \plh
\fns $\alpha_1$ and $\alpha_2$ on $S$ \st $\rho_j=\Phi^*\alpha_j$
for $j=1,2$.

In particular, if there exists a nonconstant \holo \fn with finite
energy on $X$, then there exists a surjective proper \holo mapping
$\Phi\colon X\to S$ of $X$ onto a Riemann surface~$S$ with
$\Phi_*\ol_X=\ol_S$.
\end{cor}
\begin{rmk}
Two real-valued \plh \fns $u$ and $v$ on a \con \cpx manifold have
linearly dependent differentials (i.e.~the \fns $u$,~$v$,~and~$1$
are linearly dependent) if and only if $du\wedge dv\equiv 0$.
\end{rmk}
\begin{pf*}{Proof of Corollary~\ref{L2 Castelnuovo for partial exact
forms cor}} If $\partial\rho_1$ and $\partial\rho_2$ are linearly
independent, then we may apply Theorem~\ref{L2 Castelnuovo de
Franchis thm from intro} to this pair of \holo $1$-forms. If not,
then there exist constants $\zeta_1,\zeta_2\in\C\sm\set{0}$ \st
the \fn $h=\zeta_1\rho_1+\zeta_2\rho_2\colon X\to\C$ is a
nonconstant \holo \fn with finite energy. The closed \holo
$1$-forms $\omega_1\equiv dh$ and $\omega_2\equiv hdh=2\inv
d(h^2)$ are then linearly independent and $\omega_1$ is in~$L^2$,
so we may again apply Theorem~\ref{L2 Castelnuovo de Franchis thm
from intro}. In either case, we get a proper \holo mapping
$\Phi\colon\Omega\to S$ of $X$ onto a Riemann surface~$S$ with
$\Phi_*\ol_X=\ol_S$ and the \plh \fns $\rho_1$ and $\rho_2$
descend to \plh \fns $\alpha_1$ and $\alpha_2$, respectively,
on~$S$.
\end{pf*}

\begin{defn}\label{Bdd geom weak 1-complete along set defn}
For $S\subset X$ and $k$ a positive integer, we will say that a
Hermitian manifold~$(X,g)$ {\it has bounded geometry of order $k$
along $S$} if, for some constant $C>0$ and for every point $p\in
S$, there is a biholomorphism $\Psi $ of the unit ball
$B=B(0;1)\subset \C ^n$ onto a \nbd of $p$ in $X$ \st $\Psi (0)=p$
and \stns, on $B$,
\[
C\inv g_{\C ^n}\leq\Psi^*g\leq Cg_{\C^n}
\quad\text{and}\quad|D^m\Psi ^*g|\leq C\text{ for
}m=0,1,2,\dots,k.
\]
\end{defn}

Slight modifications of the proofs of Theorem~\ref{L2 Castelnuovo
de Franchis thm from intro} and Corollary~\ref{L2 Castelnuovo for partial exact forms cor}
give the following useful generalizations:
\begin{thm}\label{L2 Castel de Franch on open set thm }
Let $\Omega$ be a nonempty domain in a \con complete Hermitian
manifold~$(X,g)$ and let $\omega_1$ and $\omega_2$ be linearly
independent closed \holo $1$-forms on~$\Omega$ \st $X$ has bounded
geometry along~$\Omega$, $g\restrict\Omega$ is K\"ahler,
$\omega_1$ is in~$L^2$, $\omega_1\wedge\omega_2\equiv 0$
on~$\Omega$, and the levels of the associated \holo mapping
\(f=({\omega_1}/{\omega_2})\colon\Omega\to\mathbb P^1\)
are closed relative to~$X$. Then there exist a surjective proper
\holo mapping $\Phi\colon\Omega\to S$ of $\Omega$ onto a Riemann
surface~$S$ with $\Phi_*\ol_\Omega=\ol_S$ and \holo $1$-forms
$\theta_1$ and $\theta_2$ on $S$ \st $\omega_j=\Phi^*\theta_j$ for
$j=1,2$.
\end{thm}
\begin{cor}\label{L2 Castelnuovo on open set for partial exact forms cor}
Let $\Omega$ be a nonempty domain in a \con complete Hermitian
manifold~$(X,g)$ and let $\rho_1$ and $\rho_2$ be two real-valued
\plh \fns on $\Omega$ \st $d\rho_1$ and $d\rho_2$ are linearly
independent, $X$ has bounded geometry along~$\Omega$,
$g\restrict\Omega$ is K\"ahler, $\rho_1$ has finite energy,
$\partial\rho_1\wedge\partial\rho_2\equiv 0$ on~$\Omega$, and the
closure (relative to~$X$) of each leaf of the (singular) \holo
foliation determined by $\partial\rho_1$ (and $\partial\rho_2$) is
contained in~$\Omega$. Then there exist a surjective proper \holo
mapping $\Phi\colon\Omega\to S$ of $\Omega$ onto a Riemann
surface~$S$ with $\Phi_*\ol_\Omega=\ol_S$ and real-valued \plh
\fns $\alpha_1$ and $\alpha_2$ on $S$ \st $\rho_j=\Phi^*\alpha_j$
for $j=1,2$.

In particular, if there exists a nonconstant \holo \fn with finite
energy on $\Omega$ whose levels are closed relative to~$X$, then
there exists a surjective proper \holo mapping
$\Phi\colon\Omega\to S$ of $\Omega$ onto a Riemann surface~$S$
with $\Phi_*\ol_\Omega=\ol_S$.
\end{cor}

\section{Proof of the cup product lemma}\label{cup
product lemma section}

Throughout this section $(X,g)$ will denote a \con complete
Hermitian manifold of dimension~$n$ with associated real
$(1,1)$-form~$\eta$. As in~\cite{Gaffney},  fixing a point $p\in
X$ and setting
\[
\tau(s)=\left\{
\begin{aligned}
1&\quad\text{if }s\leq 1\\
2-s&\quad\text{if }1<s<2\\
0&\quad\text{if }2\leq s
\end{aligned}\right.
\]
and
\[
\tau_r(x)=\tau\left(\frac{\dist(p,x)}{r}\right)
\]
for each point $x\in X$ and each number $r>0$, we get a collection
of nonnegative Lipschitz \cont \fns $\seq\tau r_{r>0}$ \stns, for
each $r>0$, we have $0\leq\tau_r\leq 1$ on $X$, $\tau_r\equiv 1$
on $B(p;r)$, $\tau_r\equiv 0$ on~$X\sm B(p;2r)$, and
$|d\tau_r|_g\leq 1/r$. Finally, for each $R>0$, $\cal M_R$ will
denote the operator given by
\[
\cal M_R(\vphi)(x)=\left\{
\begin{aligned}
\vphi(x)&\quad\text{if }|\vphi(x)|\leq R\\
R&\quad\text{if }\vphi(x)>R\\
-R&\quad\text{if }\vphi(x)<-R
\end{aligned}\right.
\]
for every (extended) real-valued \fn~$\vphi$.

\begin{pf*}{Proof of Theorem~\ref{Cup product lemma general thm
intro}} Clearly, we may assume that $n=\dim X>1$.
Assuming $X$ is K\"ahler, let $\omega_1$ and $\omega_2$ be
closed \holo $1$-forms on $X$ \st $\omega_1$ is bounded,
$\real(\omega_1)$ is exact, and $\omega_2$ is in~$L^2$. In
particular, we may fix a real-valued \plh \fn $\rho$ on $X$ \st
$\real(\omega_1)=d\rho$. Setting $d^c=-\sqrt{-1}(\partial-\dbar)$,
we get
\[
0\leq\sqrt{-1}\omega_1\wedge\overline{\omega_1}=dd^c(\rho^2)=2d(\rho
d^c\rho),
\]
and hence $\gamma=d\theta$, where $\gamma$ is the nonnegative form
of type $(n,n)$ given by
\[
\gamma\equiv\left(\sqrt{-1}\omega_1\wedge\overline{\omega_1}\right)\wedge
\left(\sqrt{-1}\omega_2\wedge\overline{\omega_2}\right)\wedge\eta^{n-2}
\]
and
\[
\theta\equiv 2\rho
(d^c\rho)\wedge(\sqrt{-1}\omega_2\wedge\overline{\omega_2})\wedge\eta^{n-2}.
\]
For every $R>0$, let $\gamma_R$ be the product of $\gamma$ and the
characteristic \fn of
\[
\setof{x\in X}{|\rho(x)|\leq R},
\]
let $\rho_R=\cal M_R(\rho)$, and let $\theta_R$ be the $L^1$
Lipschitz \cont form given by
\[
\theta_R\equiv 2\rho_R
(d^c\rho)\wedge(\sqrt{-1}\omega_2\wedge\overline{\omega_2})\wedge\eta^{n-2}.
\]
Then, for almost every $R>0$, $\gamma_R$ is equal almost
everywhere to $d\theta_R$; in fact, $\gamma_R=d\theta_R$ on
$X\sm\rho\inv(\set{\pm R})$. For each such fixed $R>0$ and each
$r>0$, Stokes' theorem gives.
\[
\int_X\tau_r\gamma_R=-\int_Xd\tau_r\wedge\theta_R.
\]
Letting $r\to\infty$ and applying the dominated convergence
theorem on the right-hand side, we get
\[
\int_X\gamma_R=0.
\]
We have $\gamma_R\geq 0$, and, therefore, $\gamma_R=0$, on
$X\sm\rho\inv(\set{\pm R})$. Letting $R\to\infty$, we get
$\gamma\equiv 0$ on $X$ and it follows that
$\omega_1\wedge\omega_2\equiv 0$.
\end{pf*}

Similar arguments yield generalizations; several examples of
which will be considered in Section~\ref{Cup product lemma general
sect}. For now, we consider two slight
generalizations of Lemma~2.7 of \cite{NR Filtered ends} which
will also allow us to give a simplified proof of the main result
of~\cite{NR Filtered ends} (see Sections~\ref{Filtered ends prelim
sect}~and~\ref{pf of filtered ends sect}).  The proof given below
is also simpler than the proof of Lemma~2.7 of \cite{NR
Filtered ends} given in that paper.

\begin{thm}\label{cup product plh fn and form on superlevel thm}
Let $\omega_1$ and $\omega_2$ be two closed \holo $1$-forms on a
domain $Y\subset X$ \st $\real(\omega_1)=d\rho_1$ for some
real-valued \plh \fn $\rho_1$ on~$Y$. Assume that, for some
constant~$a$ with $\inf\rho_1<a<\sup\rho_1$ and some \comp
$\Omega$ of $\setof{x\in Y}{a<\rho_1(x)}$, we have the following:
\begin{enumerate}
\item[(i)] $\overline\Omega\subset Y$;

\item[(ii)] The metric~$g\restrict\Omega$ is K\"ahler;

\item[(iii)] The form $\omega_1\restrict\Omega$ is bounded; and

\item[(iv)] $\int_\Omega|\omega_2|^2_g\,dV_g <\infty$.

\end{enumerate}
Then $\omega_1\wedge\omega_2\equiv 0$ on $Y$. Furthermore, if
$\omega_1$ and $\omega_2$ are linearly independent and $(X,g)$ has
bounded geometry along $\Omega$, then there exist a surjective
proper \holo mapping $\Phi\colon\Omega\to S$ of $\Omega$ onto a
Riemann surface~$S$ with $\Phi_*\ol_\Omega=\ol_S$ and \holo
$1$-forms $\theta_1$ and $\theta_2$ on $S$ \st
$\omega_j\restrict\Omega=\Phi^*\theta_j$ for $j=1,2$.
\end{thm}
\begin{pf}
Clearly, we may assume that $n=\dim X>1$.
Let $\gamma$ be the nonnegative form of type $(n,n)$ on~$Y$ given
by
\[
\gamma\equiv\left(\sqrt{-1}\omega_1\wedge\overline{\omega_1}\right)\wedge
\left(\sqrt{-1}\omega_2\wedge\overline{\omega_2}\right)\wedge\eta^{n-2}.
\]
Fixing a regular value~$b$ for $\rho_1$ with
$a<b<\sup_\Omega\rho_1$, setting
$\Omega_b=\setof{x\in\Omega}{b<\rho_1(x)}\neq\emptyset$, and
setting
\[
\theta\equiv 2(\rho_1-b)
(d^c\rho_1)\wedge(\sqrt{-1}\omega_2\wedge\overline{\omega_2})\wedge\eta^{n-2},
\]
we get $\gamma=d\theta$ on $\Omega$. For every $R>0$, let
$\gamma_R$ be the product of $\gamma$ and the characteristic \fn
of
\[
\setof{x\in Y}{|\rho_1(x)-b|\leq R},
\]
let $\alpha_R=\cal M_R(\rho_1-b)$, and let $\theta_R$ be the $L^1$
Lipschitz \cont form on~$\Omega$ given by
\[
\theta_R\equiv 2\alpha_R
(d^c\rho_1)\wedge(\sqrt{-1}\omega_2\wedge\overline{\omega_2})\wedge\eta^{n-2}.
\]
Then, for almost every $R>0$, $\gamma_R$ is equal almost
everywhere to $d\theta_R$ in~$\Omega$; in fact,
$\gamma_R=d\theta_R$ on $\Omega\sm\rho_1\inv(\set{b\pm R})$. For
each such fixed $R>0$ and each $r>0$, Stokes' theorem gives
\[
\int_{\Omega_b}\tau_r\gamma_R=-\int_{\Omega_b}d\tau_r\wedge\theta_R;
\]
since $\alpha_R\equiv 0$ on $\partial\Omega_b$. Letting
$r\to\infty$ and applying the dominated convergence theorem on the
right-hand side, we get
\[
\int_{\Omega_b}\gamma_R=0.
\]
We have $\gamma_R\geq 0$, and, therefore, $\gamma_R=0$, on
$\Omega_b\sm\rho_1\inv(\set{b\pm R})$. Letting $R\to\infty$, we get
$\gamma\equiv 0$ on $\Omega_b$ and it follows that
$\omega_1\wedge\omega_2\equiv 0$ on~$Y$.

Assume now that $\omega_1$ and $\omega_2$ are linearly independent
and $(X,g)$ has bounded geometry along $\Omega$. Since $\rho_1$ is
constant on the levels of the \holo map \(f=\omega_1/\omega_2\),
those levels which meet~$\Omega$ are contained in~$\Omega$. Thus
Theorem~\ref{L2 Castel de Franch on open set thm } gives the
desired proper \holo mapping to a Riemann surface.
\end{pf}
Applying the above theorem together with Corollary~\ref{L2
Castelnuovo on open set for partial exact forms cor}, we get the
following:
\begin{cor}\label{cup product two plh fns on superlevel cor}
Let $\rho_1$ and $\rho_2$ be two real-valued \plh \fns on a domain
$Y\subset X$. Assume that, for some constant~$a$ with
$\inf\rho_1<a<\sup\rho_1$ and some \comp $\Omega$ of $\setof{x\in
Y}{a<\rho_1(x)}$, we have the following:
\begin{enumerate}
\item[(i)] $\overline\Omega\subset Y$,

\item[(ii)] The metric~$g\restrict\Omega$ is K\"ahler,

\item[(iii)] The form $d\rho_1\restrict\Omega$ is bounded, and

\item[(iv)] $\int_\Omega|d\rho_2|^2_g\,dV_g <\infty$.

\end{enumerate}
Then $\partial\rho_1\wedge\partial\rho_2\equiv 0$ on $Y$.
Furthermore, if $d\rho_1$ and $d\rho_2$ are linearly independent
and $(X,g)$ has bounded geometry along $\Omega$, then there exist
a surjective proper \holo mapping $\Phi\colon\Omega\to S$ of
$\Omega$ onto a Riemann surface~$S$ with $\Phi_*\ol_\Omega=\ol_S$
and \plh \fns $\alpha_1$ and $\alpha_2$ on $S$ \st
$\rho_j\restrict\Omega=\Phi^*\alpha_j$ for $j=1,2$.
\end{cor}

\section{An application to filtered ends of K\"ahler manifolds}\label{Filtered ends prelim sect}

Let $X$ be a \con complete K\"ahler manifold.  According to
\cite{Gro-Sur la groupe fond}, \cite{Li Structure complete
Kahler}, \cite{Gro-Kahler hyperbolicity}, and Theorem~3.4 of
\cite{NR-Structure theorems}, if $X$ has at least $3$~ends, and
either $X$ has bounded geometry of order~$2$ or $X$ is weakly
$1$-complete or $X$ admits a positive symmetric Green's \fn which
vanishes at infinity, then $X$ maps properly and holomorphically
onto a Riemann surface. The ends condition was weakened in
\cite{Delzant-Gromov Cuts} and \cite{NR Filtered ends} to the
condition that $X$ have at least $3$~\textit{filtered ends}
relative to the universal covering.  The techniques and results described in
the previous sections allow one to simplify the proof of the main result
of~\cite{NR Filtered ends} (see Theorem~0.1 and
Theorem~3.1 of~\cite{NR Filtered ends}) in the following sense. The proof given
in~\cite{NR Filtered ends} relied heavily on a weak version of
Theorem~\ref{cup product plh fn and form on superlevel thm} in which both of the
\holo $1$-forms are assumed to be in~$L^2$ on the domain~$\Omega$ and to have
exact real parts (see Lemma~2.7 of~\cite{NR Filtered ends}). The proof
of Theorem~\ref{cup product plh fn and form on superlevel thm} given in
Section~\ref{cup product lemma section} is simpler than that of the weak
version given in~\cite{NR Filtered ends}. Moreover,
Theorem~\ref{cup product plh fn and form on superlevel thm}, being stronger, allows one to eliminate
some of the technical arguments used in~\cite{NR Filtered ends}. In fact, one can avoid any
direct use of the general theory of massive sets due to Grigor'yan~\cite{Grigoryan};
a central technique employed in~\cite{NR Filtered ends}.
In this section, we recall the
required definitions and preliminary facts. The new proof appears
in Section~\ref{pf of filtered ends sect}.

\begin{defn}\label{ends filtered ends def}
Let $M$ be a \con manifold.
\begin{enumerate}
\item[(a)] By an {\it end} of~$M$, we will mean either a \comp $E$
of $M\setminus K$ with non\cpt closure, where $K$ is a given \cpt
subset of $M$, or an element of
$$
\lim _{\leftarrow } \pi _0 (M\setminus K),
$$
where the limit is taken as $K$ ranges over the \cpt subsets of
$M$ (or the \cpt subsets of $M$ whose complement $M\setminus K$
has no \rel \cpt \compsns). The number of ends of $M$ will be
denoted by~$e(M)$. For a \cpt set $K$ \st $M\setminus K$ has no
\rel \cpt \compsns, we will call
$$
M\setminus K=E_1\cup \cdots \cup E_m,
$$
where $E_1, \dots , E_m$ are the distinct \comps of $M\setminus
K$, an {\it ends decomposition} for~$M$.

\item[(b)] (Following Geoghegan \cite{Geoghegan}) For $\Upsilon
:\widetilde M\to M$ the universal covering of $M$, elements of the
set
$$
\lim_{\leftarrow }\pi_0 [\Upsilon\inv(M\setminus K)],
$$
where the limit is taken as $K$ ranges over the \cpt subsets of
$M$ (or the \cpt subsets of $M$ whose complement $M\setminus K$
has no \rel \cpt \compsns) will be called {\it filtered ends}. The
number of filtered ends of $M$ will be denoted by~$\tilde e(M)$.
\end{enumerate}
\end{defn}

Clearly, $\tilde e(M)\geq e(M)$. In fact, for $k\in \N$, we have
$\tilde e(M)\geq k$ if and only if there exists an ends
decomposition $M\setminus K=E_1\cup \cdots \cup E_m$ for $M$
\stns, for $\Gamma_j=\image{\pi_1(E_j)\to\pi_1(M)}$ for
$j=1,\dots,m$, we have
$$
\sum _{j=1}^m[\pi _1(M):\Gamma _j]\geq k.
$$
Moreover, if $\widehat M\to M$ is a \con covering space, then
$\tilde e(\widehat M)\leq \tilde e(M)$ with equality if the
covering is finite.

\begin{defn}\label{weakly 1-complete along set def}
We will say that a \cpx manifold $X$ is {\it weakly $1$-complete
along} a subset~$S$ if there exists a \cont \plsh \fn $\vphi$ on
$X$ \st
$$
\setof{x\in S}{\vphi (x)<a}\Subset X \quad \forall \, a\in \R .
$$
\end{defn}

\begin{defn}\label{Special end definition}
We will call an end $E$ of a \con non\cpt complete
Hermitian manifold $(X,g)$ {\it special} if $E$ is of at least one
of the following types:
\begin{enumerate}
\item[(BG)] $(X,g)$ has bounded geometry of order $2$ along $E$;

\item[(W)] $X$ is weakly $1$-complete along $E$;

\item[(RH)] $E$ is a hyperbolic end and the Green's \fn vanishes
at infinity along $E$; or

\item[(SP)] $E$ is a parabolic end, the Ricci curvature of $g$ is
bounded below on $E$, and there exist positive constants $R$ and
$\delta $ \st
$$
\Vol \big( B(p;R)\big) >\delta \quad \forall \, p\in E.
$$
\end{enumerate}
An ends decomposition for $X$ in which each of the ends is special
will be called a {\it special ends decomposition}.
\end{defn}
\begin{rmks}
1. (BG) stands for ``bounded geometry,'' (W) for ``weakly
$1$-complete,'' (RH) for ``regular hyperbolic,'' and (SP) for
``special parabolic.''

\noindent 2. A parabolic end of type (BG) is also of type (SP).

\noindent 3. If $E$ and $E'$ are ends with $E'\subset E$ and $E$
is special, then $E'$ is special.

\noindent 4. We recall that an end $E$ of a \con Riemannian
manifold~$(M,g)$ is hyperbolic if and only if there exists a
bounded nonnegative \cont sub\harm \fn $\alpha$ on~$M$ \st
$\alpha\equiv 0$ on $M\sm E$ and $\sup_E\alpha>0$. Such a
\fnns~$\alpha$ is called an \textit{admissible sub\harm \fnns}
for~$E$ in~$M$. The end~$E$ is special of type~(RH) if and only if
we may choose~$\alpha$ so that $\alpha\to\sup\alpha$ at infinity
in~$\overline E$. As in the work of Grigor'yan~\cite{Grigoryan},
any open set~$E$ (whether or not it's an end) is called
\textit{massive} if there exists an admissible sub\harm
\fn for~$E$. General massive sets are applied in~\cite{NR Filtered
ends}, but the results of Section~\ref{cup product lemma section}
will allow us to restrict our attention to hyperbolic ends.
\end{rmks}

Special ends in a complete K\"ahler manifold allow one to produce
\plh \fns and, in some cases, \holo \fnsns.  In particular, one
gets the following:
\begin{thm}[\cite{Gro-Sur la groupe fond}, \cite{Li Structure complete Kahler},
\cite{Gro-Kahler hyperbolicity}, and Theorem~3.4 of
\cite{NR-Structure theorems}]\label{Old 3 ends theorem} If $(X,g)$
is a \con complete K\"ahler manifold which admits a special ends
decomposition and $e(X)\geq 3$, then $X$ admits a proper \holo
mapping onto a Riemann surface.
\end{thm}

The main result of~\cite{NR Filtered ends} is the following
generalization (see Theorem~3.1 of~\cite{NR Filtered ends}):
\begin{thm}\label{filtered ends theorem} If $(X,g)$
is a \con complete K\"ahler manifold which admits a special ends
decomposition and $\til e(X)\geq 3$, then $X$ admits a proper
\holo mapping onto a Riemann surface.
\end{thm}

The goal of this section and Section~\ref{pf of filtered ends
sect} is to describe a simpler proof of the above fact. We will
produce independent \plh \fns by applying Theorem~2.6
of~\cite{NR-Structure theorems}, which is contained implicitly in
the work of Sario, Nakai, and their collaborators
\cite{Nakai1},\cite{Nakai2}, \cite{SaNa}, \cite{SaNo}, \cite{RS}
and the work of Sullivan~\cite{Sul} (see also \cite{Li Structure
complete Kahler} and \cite{LT}). This fact is also applied
in~\cite{NR Filtered ends} along with the more general theory of
massive sets~\cite{Grigoryan}, but we will not need general
massive sets in this paper. In fact, we will only need the
following weak version of Theorem~2.6 of~\cite{NR-Structure
theorems}:

\begin{thm}\label{plh fn along ends modified thm}
Let $(X,g)$ be a \con complete K\"ahler manifold with
an ends decomposition $X\sm K=E_1\cup\cdots\cup E_m$ \st $m>1$ and
\stns, for each $j=1,\dots,m$, $E_j$ is a hyperbolic end or a
special end of type~(SP). Then there exists a \plh \fn $\rho\colon
X\to\R$ \stns, for each $j=1,\dots,m$, we have the following:
\begin{enumerate}
\item[(i)] If $E_j$ is a hyperbolic end, then
$0<\rho\restrict{E_j}<1$ and $\rho\restrict{E_j}$ has finite
energy;

\item[(ii)] If $E_1$ is a hyperbolic end (a special end of
type~(RH)), then
\[
\limsup_{x\to\infty}\rho\restrict{\overline{E_1}}(x)=1\qquad\text{(respectively,
}\lim_{x\to\infty}\rho\restrict{\overline {E_1}}(x)=1\text{)};
\]
and

\item[(iii)] If $E_1$ is a special end of type~(SP), then
\[
\lim_{x\to\infty}\rho\restrict{\overline{E_1}}(x)=\infty.
\]

\end{enumerate}
\end{thm}
\begin{rmk}
Theorem~2.6 of \cite{NR-Structure theorems} is stated for
dimension~$n>1$, but it actually holds in arbitrary dimension.
On the other hand, we will only need Theorem~\ref{plh fn along
ends modified thm} for $n>1$.
\end{rmk}
\begin{pf*}{Proof of Theorem~\ref{plh fn along ends modified thm}}
Applying Theorem~2.6 of~\cite{NR-Structure theorems}, we get a
nonconstant \plh \fn $\alpha\colon X\to\R$ \stns, for each
$j=1,\dots,m$, we have the following:
\begin{enumerate}
\item[(\ref{plh fn along ends modified thm}.1)] If $E_j$ is a
hyperbolic end, then $\alpha\restrict{E_j}$ is bounded with finite
energy and
\[
\liminf_{x\to\infty}\alpha\restrict{\overline{E_j}}(x)\left\{
\begin{aligned}
=0&\qquad\text{if }j=1\\
>0&\qquad\text{if }j>1
\end{aligned}
\right.
\]

\item[(\ref{plh fn along ends modified thm}.2)] If $E_j$ is a
special end of type~(RH), then
\[
\lim_{x\to\infty}\alpha\restrict{\overline{E_j}}(x)=\left\{
\begin{aligned}
0&\qquad\text{if }j=1\\
1&\qquad\text{if }j>1
\end{aligned}
\right.
\]

\item[(\ref{plh fn along ends modified thm}.3)] If $E_j$ is a
special end of type~(SP), then
\[
\lim_{x\to\infty}\alpha\restrict{\overline{E_j}}(x) =\left\{
\begin{aligned}
\infty&\qquad\text{if }j=1\\
\infty&\qquad\text{if }j>1\text{ and }X\text{ is hyperbolic (i.e.
}E_i\text{ is hyperbolic for some }i\text{)}\\
-\infty&\qquad\text{if }j>1\text{ and }X\text{ is parabolic (i.e.
}E_1,\dots,E_m\text{ are parabolic)}
\end{aligned}
\right.
\]

\end{enumerate}
Let $H$ be the union of all of those ends $E_j$ which
are hyperbolic and fix $s\in\R$ with $\alpha<s$ on~$H$. If $H\neq\emptyset$, then
the maximum principle implies that $\alpha>0$ on~$X$. Thus the
\fn
\[
\rho\equiv\left\{
\begin{aligned}
1-(\alpha/s)&\qquad\text{if }H\neq\emptyset\text{ and }E_1\text{
is hyperbolic}\\
\alpha/s&\qquad\text{if }H\neq\emptyset\text{ and
}E_1\text{ is
parabolic}\\
\alpha&\qquad\text{if }H=\emptyset
\end{aligned}
\right.
\]
has the required properties.
\end{pf*}

The following lemma may be viewed as a consequence of Theorem~\ref{Old 3 ends theorem}
(see, for example, the proof of
Theorem~4.6 of \cite{NR-Structure theorems}):
\begin{lem}\label{Open set to Riemann surface gives global lemma}
Let $(X,g)$ be a \con complete K\"ahler manifold which is \cpt or
which admits a special ends decomposition. If some nonempty open
subset of $X$ admits a surjective proper \holo mapping onto a
Riemann surface, then $X$ admits a surjective proper \holo mapping
onto a Riemann surface.
\end{lem}

The following easy observation will enable us to produce \plh \fns
by passing to a covering.
\begin{lem}\label{Lifting hyperbolic to cover} Let $(X,g)$ be a \con
complete K\"ahler manifold, let $\Upsilon\colon\what X\to X$ be a
\con covering space, and let $\hat g=\Upsilon^*g$, .
\begin{enumerate}

\item[(a)] If $E_1$ is a hyperbolic end of $X$, then any end $F$
of $\what X$ containing a \compns~$E$ of $\Upsilon\inv(E_1)$ is a
hyperbolic end.

\item[(b)] If $E_1$ is a special end with smooth boundary and $E$
is a \comp of $\Upsilon\inv(E_1)$ for which the restriction $E\to
E_1$ is a finite covering, then $E$ is a special end of the same
type.

\item[(c)] If $X\sm K=E_1\cup\cdots\cup E_m$ is an ends
decomposition into hyperbolic ends with smooth boundary and $E$ is
a \comp of $\Upsilon\inv(E_1)$ for which the restriction $E\to
E_1$ is a finite covering, then every \comp of $\what X\sm\partial
E$ with non\cpt closure is a hyperbolic end of~$\what X$.

\end{enumerate}
\end{lem}
\begin{pf} Let $X\sm K=E_1\cup\cdots\cup E_m$ be an ends
decomposition, let $\what E_j=\Upsilon\inv(E_j)$ for
$j=1,\dots,m$, and let $E$ be a \comp of $\what E_1$. If $E_1$ is
a hyperbolic end, $\alpha$ is an admissible sub\harm \fn
for~$E_1$, and $F$ is an end of $\what X$ containing~$E$, then we
have the admissible sub\harm \fn
\[
\beta\equiv\left\{
\begin{aligned}
\alpha\circ\Upsilon&\qquad\text{on }E\\
0&\qquad\text{on }\what X\sm E
\end{aligned}
\right.
\]
for~$F$. Thus (a) is proved.

If $E_1$ is a smooth domain and the restriction $E\to E_1$ is a
finite covering, then $E$ is an end of~$\what X$ and there exist \nbds $V$ and $V_1$ of
$\overline E$ and $\overline{E_1}$, respectively, \st
$V\cap\Upsilon\inv(E_1)=E$ and $V\to V_1$ is also a finite
covering space. Clearly, if $E_1$ is special of type~(BG),~(W),~or~(RH),
then $E$ is special of the same type. If $E$ is a hyperbolic end with
an admissible sub\harm \fnns~$\alpha$, then the \fn
\[
\alpha_1(x)\equiv\left\{
\begin{aligned}
\sum_{y\in\Upsilon\inv(x)\cap E}\alpha(y)&\qquad\text{if }x\in E_1\\
0&\qquad\text{if }x\in X\sm E_1
\end{aligned}
\right.
\]
is an admissible sub\harm \fn for~$E_1$. It now follows easily that,
if $E_1$ is special of type~(SP), then~$E$ must also be special of type~(SP).
Thus~(b) is proved.

Finally, suppose that $E_j$ is a hyperbolic end and a smooth
domain for each $j=1,\dots,m$ and that the restriction $E\to E_1$
is a finite covering. In particular, forming $V\to V_1$ as above,
we see that $\partial E$ is \cpt and every \compns~$F$ of $\what
X\sm\partial E$ with non\cpt closure is an end. Furthermore, $F$
must meet $\what E_j$ for some~$j$. For if not, then $F$ must be a
\concomp of $\Upsilon\inv(X\sm\overline{E_1})$ contained in $\what
K=\Upsilon\inv(K)$. Thus we have the \con covering space
$\overline F\to \overline{F_1}$ and the covering space $\partial F\to\partial F_1$
for some \compns~$F_1$ of~$X\sm\overline{E_1}$ contained in~$K$
(here, we have used smoothness). Since $\partial
F\subset\partial E$ and $\partial E\to\partial E_1$ is a finite
covering space (of manifolds), we see that $\overline
F\to\overline{F_1}\subset K$ is a finite cover, which contradicts
the noncompactness of~$\overline F$. Thus $F$ must meet and,
therefore, contain, $\what E_j$ for some~$j$ and~(a)
then implies that $F$ is a hyperbolic end.  Thus~(c) is proved.
\end{pf}

\section{Proof of the filtered ends result}\label{pf of filtered
ends sect}

The first step in the new proof of Theorem~\ref{filtered ends
theorem} is to reduce to the case in which all of the ends of the
manifold are hyperbolic special ends of type~(BG). Toward this
goal, we first recall the following two facts:

\begin{lem}[See Lemma~3.2 of \cite{NR Filtered ends}]\label{filtered ends of subdomains lemma}
Let $M$ be a \con non\cpt $\cinf$ manifold and let $k\in\N$.
\begin{enumerate}
\item[(a)] Given an end $E$ in $M$ with
$\left[\pi_1(M):\image{\pi_1(E)\to\pi_1(M)}\right]\geq k$, there
exists a \cpt set $D\subset M$ \stns, if $\Omega$ is a domain
containing $D$, then $\Omega\cap E$ is an end of $\Omega$ and, for
any end $F$ of~$\Omega$ contained in $E$, we have
$\left[\pi_1(\Omega):\image{\pi_1(F)\to\pi_1(\Omega)}\right]\geq
k$.

\item[(b)] If $\til e(M)\geq k$, then there exists a \cpt set
$D\subset M$ \stns, for every domain $\Omega$ containing $D$, we
have $\til e(\Omega)\geq k$.

\end{enumerate}
\end{lem}

\begin{lem}[See Lemma~3.3 of \cite{NR Filtered ends}]\label{Type W to RH lemma}
Let $(X,g)$ be a \con complete K\"ahler manifold, let $E$ be a
special end of type~(W) in $X$, let $k,l\in\N$ with
\[
\til e(X)\geq k\qquad\text{and}\qquad
\left[\pi_1(X):\image{\pi_1(E)\to\pi_1(X)}\right]\geq l,
\]
and let $D$ be a \cpt subset of $X$. Then there exists a
domain $X'$ in $X$, a complete K\"ahler metric $g'$ on $X'$, a
\cpt set $K\subset X'$, and disjoint domains $E_0,\dots,E_m$ \st
\begin{enumerate}

\item[(i)] $(X\sm E)\cup D\subset E_0$, $X'\sm K=E_0\cup
E_1\cup E_2\cup\cdots\cup E_m$, and $E_0\cap E\Subset X'$;

\item[(ii)] On $E_0$, $g'=g$;

\item[(iii)]  For each $j=1,\dots,m$, $E_j$ is a special end of
type~(RH) and (W) satisfying $\left[\pi_1(X'):\text{\rm
im}\,[\pi_1(E_j)\to\pi_1(X')]\right]\geq l$; and

\item[(iv)] $\til e(X')\geq k$.

\end{enumerate}
\end{lem}

The above lemmas allow us to replace special ends of type~(W) with
special ends of type~(RH). The following lemma will allow us to
replace special ends of type~(RH) with special ends of type~(BG)
(under the right conditions).

\begin{lem}\label{Type RH to BG lemma}
Let $(X,g)$ be a \con complete K\"ahler manifold, let $E$ be an
end of $X$, let $k,l\in\N$ with
\[
\til e(X)\geq k\qquad\text{and}\qquad
\left[\pi_1(X):\image{\pi_1(E)\to\pi_1(X)}\right]\geq l,
\]
and let $D$ be a \cpt subset of $X$. Assume that, for some
$R\in(0,\infty]$, there exists a \cont \fn $\rho\colon\overline
E\to(0,R)$ \st $\rho$ is \plh on~$E$ and
\[
\lim_{x\to\infty}\rho\restrict{\overline E}(x)=R
\]
(in particular, $E$ is a special end of type~(W)). Then there
exists a domain $X'$ in $X$, a complete K\"ahler metric $g'$ on
$X'$, a \cpt set $K\subset X'$, and disjoint domains
$E_0,\dots,E_m$ \st
\begin{enumerate}

\item[(i)] $(X\sm E)\cup D\subset E_0$, $X'\sm K=E_0\cup
E_1\cup E_2\cup\cdots\cup E_m$, and $E_0\cap E\Subset X'$;

\item[(ii)] We have $g'\geq g$ on $X'$ and $g'=g$ on~$E_0$;

\item[(iii)]  For each $j=1,\dots,m$, $E_j$ is a special end of
type~(BG), (RH), and (W) for~$(X',g')$ satisfying
$\left[\pi_1(X'):\text{\rm
im}\,[\pi_1(E_j)\to\pi_1(X')]\right]\geq l$; and

\item[(iv)] $\til e(X')\geq k$.

\end{enumerate}
\end{lem}
\begin{pf}
By Lemma~\ref{filtered ends of subdomains lemma}, we may assume
without loss of generality that $D$ is nonempty and \conns; $\partial E\subset
D$; and, if $\Omega$ is any domain in $X$ containing $D$, then
$\til e(\Omega)\geq k$, $\Omega\cap E$ is an end of $\Omega$, and,
for any end $F$ of~$\Omega$ contained in $E$, we have
$\left[\pi_1(\Omega):\image{\pi_1(F)\to\pi_1(\Omega)}\right]\geq
l$. Fixing positive constants $a$,~$b$,~and~$c$ with
$\max_{D\cap\overline E}\rho<a<b<c<R$ and a \cinf \fn
$\chi\colon\R\to\R$ \st $\chi'\geq 0$ and $\chi''\geq 0$ on~$\R$,
$\chi(t)=0$ for $t\leq a$, and $\chi(t)=t-b$ for $t\geq c$, we get
a \cinf \plsh \fn
\[
\vphi\equiv\left\{
\begin{aligned}
\chi(\rho)&\qquad\text{on }E\\
0&\qquad\text{on }X\sm E
\end{aligned}
\right.
\]
\st $0\leq\vphi<R-b$, $\vphi\equiv 0$ on a \nbd of $(X\sm E)\cup
D$, and $\vphi=\rho-b$ on the complement in~$E$ of the \cpt set
$\setof{x\in \overline E}{\vphi(x)\leq c-b}$. Finally, we may fix
a regular value~$r$ for~$\vphi$ with $c-b<r<R-b$.

On the \comp $X'$ of $\setof{x\in E}{\vphi(x)<r}\cup (X\setminus
E)$ containing the \con set $(X\setminus E)\cup D$, we may form
the complete K\"ahler metric
\[
g'\equiv g+\lev{-\log(r-\vphi)}.
\]
We have $X'\cap E\Subset X$ since $\vphi\to R-b$ at infinity
in~$\overline E$. We also have $g'=g$ on the interior $V$ of
$\setof{x\in X'}{\vphi(x)=0}$ and, since $\vphi\to r$ at $\partial
X'$, the closure of the \comp $E_0$ of~$V$ containing $(X\setminus
E)\cup D$ is contained in $X'$ and the set
\[
K\equiv X'\sm\left[E_0\cup\setof{x\in X'}{\vphi(x)>0}\right]
\]
is \cptns. Furthermore, by the maximum principle, each of the
\comps $E_1,\dots,E_m$ of the nonempty set $\setof{x\in
X'}{\vphi(x)>0}$ is not \rel \cpt in $X'$ and is, therefore, a
special end of type~(RH)~and~(W) with respect to~$g'$ (in fact,
with respect to any complete K\"ahler metric on~$X'$). Finally, we
have
\[
F\equiv E_1\cup\cdots\cup E_m\Subset E
\]
and, near each point in $\partial F\cap\partial X'$, the local
defining \fn $\vphi-r$ for~$F$ in~$X$ is \plh with nonvanishing
differential.  The argument on~p.~831 in~\cite{NR-Structure
theorems} now shows that $g'$ has bounded geometry of order~$2$
(in fact, of all orders) along~$F$.
\end{pf}
\begin{rmk}
Clearly, either $X'\sm K=E_0\cup E_1\cup E_2\cup\cdots\cup E_m$ is
an ends decomposition or $E_0\Subset X'$ and $X'\sm(K\cup E_0)=
E_1\cup E_2\cup\cdots\cup E_m$ is an ends decomposition.
\end{rmk}

We may now reduce to the bounded geometry case.
\begin{lem}\label{Reduction to BG lemma}
To obtain Theorem~\ref{filtered ends theorem}, it suffices to
prove the theorem for every \con complete K\"ahler manifold~$X$
which has at least~$3$ filtered ends and which admits a special ends
decomposition into at least $\min(e(X),2)$ ends, each of which is
hyperbolic of type~(BG).
\end{lem}
\begin{pf}
Given $(X,g)$ as in the statement of Theorem~\ref{filtered ends
theorem}, we may choose a special ends decomposition $X\sm
K=E_1\cup\cdots\cup E_m$ for~$X$ \st $m\geq\min(e(X),2)$ and \stns,
setting $\Gamma_j=\image{\pi_1(E_j)\to\pi_1(X)}$ for
$j=1,\dots,m$, we have
\[
\sum_{j=1}^m[\pi_1(X):\Gamma_j]\geq 3.
\]
According to Lemma~\ref{Open set to Riemann surface gives global
lemma}, in order to obtain a proper \holo mapping of $X$ onto a
Riemann surface, it suffices to find such a mapping for some
nonempty open subset of~$X$. Applying Lemma~\ref{Type W to RH
lemma} to each end of type~(W) and working on a suitable subdomain
in place of~$X$, we see that we may assume without loss of
generality that each of the ends is of type (BG), (RH), or (SP).

Note that the parabolic ends of type~(BG) are also of type~(SP).
Thus, if $m\geq 2$, then, for each~$j$, Theorem~\ref{plh fn along
ends modified thm} provides a \plh \fnns~$\alpha$ on $X$ \st
$\alpha\restrict{\overline E_j}$ is an \exh \fn if the end is of
type~(SP) and a bounded \exh \fn if the end is of type~(RH). Thus,
for $m\geq 2$, Lemma~\ref{Type RH to BG lemma} implies that we may
assume that each end is hyperbolic of type~(BG) (the condition
$m\geq 2=\min(e(X),2)$ will be satisfied automatically for the
associated subdomain $X'$ which replaces~$X$).

Thus it suffices to consider the case in which $e(X)=1$ and $E_1$
is a special end of type~(RH) or~(SP) (i.e. $X$ is itself a
special end of type~(RH) or~(SP)). We may also choose the end
$E_1$ to be a \cinf domain. For a point $x_0\in E_1$,
$\Gamma_1\equiv\image{\pi _1(E_1,x_0)\to\pi_1(X,x_0)}$ is of index
$\geq 3$.  Thus we may fix a \con covering space $\Upsilon
:\widehat X\to X$ and a point $y_0\in \widehat X$ \st
$\Upsilon_*\pi _1(\widehat X,y_0)=\Gamma_1$. Hence $\Upsilon $
maps a \nbd of the closure of the \comp $E$ of $\widehat
E_1=\Upsilon\inv (E_1)$ containing $y_0$ isomorphically onto a
\nbd of $\overline{E_1}$ and, since $\#\Upsilon\inv (x_0)=[\pi
_1(X,x_0):\Gamma_1]\geq 3$, we have $\widehat E_1\setminus
E\neq\emptyset$.

In particular, $e(\what X)>1$. By Lemma~\ref{Lifting hyperbolic to
cover}, $E$ is a special end of the same type as~$E_1$ (type (RH)
or type~(SP)) and any \compns~$F$ of $\what X\sm\overline E$ with
non\cpt closure is either a hyperbolic end (this is the case if,
for example, if $X$ is of type~(RH)) or a special end of type (SP)
(which is the case if $X$ is of type~(SP) and $F$ is not
hyperbolic). Therefore, by Theorem~\ref{plh fn along ends modified
thm}, there exists a real-valued \plh \fn $\beta$ on~$\what X$
whose restriction to $\overline E$ is either an \exh \fn or a
bounded \exh \fnns. Applying Lemma~\ref{Type RH to BG lemma} to
the \fn
\[
\rho=\beta\circ\left(\Upsilon\restrict{\overline
E}\right)\inv-\inf_E\beta
\]
on~$\overline{E_1}$, we get the associated subdomain and K\"ahler
metric $(X',g')$ in which \textit{any} end is hyperbolic and
special of type~(BG), (RH), and (W). Choosing an ends
decomposition into at least $\min(e(X'),2)$ ends, we see that the
claim follows.
\end{pf}

The next two lemmas give the theorem when the manifold has only
finitely many filtered ends.

\begin{lem}\label{finite filtered ends for mfld lemma}
Suppose $M$ is a \con non\cpt manifold and $M\sm K=E_1\cup\cdots\cup E_m$ is an ends
decomposition \stns, setting
$\Gamma_j=\image{\pi_1(E_j)\to\pi_1(M)}$ for $j=1,\dots,m$, we
have
\[
k=\sum_{j=1}^m[\pi_1(M):\Gamma_j]<\infty.
\]
Then there exists a \con finite covering space $\what M\to M$ with
$e(\what M)\geq k$.
\end{lem}
\begin{pf}
The lifting of $M\setminus K$ to the universal covering $\wtil
M\to M$ has~$k$ \compsns, and the action of $\pi_1(M)$ permutes
these \compsns. Thus we get a homomorphism of $\pi_1(M)$ into the
symmetric group on $k$ objects, and hence the kernel $\Gamma$ is a
normal subgroup of finite index. The quotient $\widehat
M=\Gamma\big\backslash\wtil M\to M$ is, therefore, a finite
covering with at least~$k$ ends.
\end{pf}

\begin{lem}\label{finite filtered ends case of thm lemma}
Let $(X,g)$ be a \con complete K\"ahler manifold which admits a
special ends decomposition $X\setminus K=E_1\cup\cdots\cup E_m$
\stns, setting $\Gamma_j=\image{\pi_1(E_j)\to\pi_1(X)}$ for
$j=1,\dots,m$, we have
\[
3\leq\sum_{j=1}^m[\pi_1(X):\Gamma_j]<\infty.
\]
Then $X$ admits a proper \holo mapping onto a Riemann surface.
\end{lem}
\begin{pf}
By Lemma~\ref{finite filtered ends for mfld lemma}, $X$ admits a
\con finite covering space $\what X\to X$ with at least three
ends. Theorem~\ref{Old 3 ends theorem} and Lemma~\ref{Lifting hyperbolic to cover}
imply that $\what X$
admits a proper \holo mapping onto a Riemann surface. The claim
follows since any normal \cpx space which is the image of a
\holoconvex \cpx space under a proper \holo mapping is itself
\holoconvexns.
\end{pf}

\begin{pf*}{Proof of Theorem~\ref{filtered ends theorem}}
By Theorem~\ref{Old 3 ends theorem} and Lemma~\ref{Reduction to BG
lemma}, it suffices to consider a \con complete K\"ahler
manifold~$(X,g)$ \st $\til e(X)\geq 3$, $m=e(X)\leq 2$, and we
have an ends decomposition $X\sm K=E_1\cup\cdots\cup E_m$ in which
each of the ends is hyperbolic and special of type~(BG). In
particular, since $m=e(X)$, \textit{every} end in $X$ is
hyperbolic and special of type~(BG). Thus we may also choose the
ends decomposition so that each of the ends~$E_j$ is a smooth
domain and, by Lemma~\ref{finite filtered ends case of thm lemma},
we may assume that the subgroup $\Gamma\equiv\image{\pi
_1(E_1)\to\pi_1(X)}$ is of infinite index in~$\pi_1(X)$. According
to Lemma~\ref{Open set to Riemann surface gives global lemma}, it
suffices to show that some nonempty open subset of $X$ admits a
proper \holo mapping onto a Riemann surface.

Let $\Upsilon\colon\what X\to X$ be a \con covering space with
$\Upsilon_*\pi_1(\what X)=\Gamma$ and let $\hat g=\Upsilon^*g$.
Then $\Upsilon$ is an infinite covering map, $\Upsilon$ maps a
\nbd of the closure of some \compns~$E$ of $\what
E_1=\Upsilon\inv(E_1)$ isomorphically onto a \nbd
of~$\overline{E_1}$, $(\what X,\hat g)$ has bounded geometry, and
each \comp of $X\sm\partial E$ with non\cpt closure is a
hyperbolic end (by Lemma~\ref{Lifting hyperbolic to cover}). In
particular, it suffices to show that some nonempty open subset of
$\what X$ admits a proper \holo mapping onto a Riemann surface.
For then $\what X$ and, therefore, some nonempty open subset of
$E\cong E_1\subset X$, admits such a mapping.

If the restriction of $\Upsilon$ to each \comp of $\what E_1$
gives a finite covering space of $E$, then $e(\what X)=\infty$ and
Theorem~\ref{Old 3 ends theorem} gives the desired proper \holo
mapping to a Riemann surface. Thus we may assume that there is a
\comp~$Y$ of $\what E_1$ \st the restriction $Y\to E_1$ is an
infinite covering. Since each \comp of $X\sm\partial E$ with
non\cpt closure is a hyperbolic end, Theorem~\ref{plh fn along
ends modified thm} provides a finite energy \plh \fn
$\rho\colon\what X\to(0,1)$ \st
\[
\limsup_{x\to\infty}\rho\restrict{\overline E}(x)=1.
\]
The $L^\infty/L^2$-comparison for local \holo \fns on a bounded geometry
Hermitian manifold implies that~$d\rho$ is bounded. Thus
the \plh \fn
\[
\rho_1\equiv\rho\circ\left(\Upsilon\restrict E \right)\inv
\circ\Upsilon\restrict Y
\]
on~$Y$ also has bounded differential. On the other hand, $d\rho_1$
is \textit{not} in~$L^2$ because the covering space $Y\to E_1$ is
infinite, so the \holo $1$-form $\omega_1\equiv\partial\rho_1$ and
the $L^2$ \holo $1$-form $\omega_2\equiv\partial\rho\restrict Y$
are linearly independent. Fixing $a$ with
\[
\max_{\partial E}\rho<a<1,
\]
we get a nonempty \concomp $\Omega$ of $\setof{x\in
Y}{a<\rho_1(x)}$ with $\overline\Omega\subset Y$. Theorem~\ref{cup
product plh fn and form on superlevel thm} now implies that
$\Omega$ admits a proper \holo mapping onto a Riemann surface and
the theorem follows.
\end{pf*}

\section{Further generalizations of the cup product lemma}\label{Cup
product lemma general sect}

In this section, the techniques described in Section~\ref{cup
product lemma section} will be extended to give several different
versions of the cup product lemma. Throughout this section,
$(X,g)$ will denote a \con complete Hermitian manifold of
dimension~$n>1$ with associated real $(1,1)$-form~$\eta$. We may
also form the collection of \fns $\seq\tau r_{r>0}$ and the
operator $\cal M_R$ for $R>0$ as in the beginning of
Section~\ref{cup product lemma section}.

Suppose $(\cal V,g)$ is a Hermitian inner product space, $\eta$ is
the skew-symmetric real form of type~$(1,1)$ associated to~$g$,
and $A$ is a Hermitian symmetric form on~$\cal V$ with associated
skew-symmetric real $(1,1)$-form~$\alpha$. We will write $A\geq 0$ and
$\alpha\geq 0$ ($A>0$ and $\alpha>0$) if $A$ is nonnegative
definite (respectively, positive definite). Given a positive
integer $q$, we will write $A\geqtrace gq 0$ and $\alpha\geqtrace
gq 0$ ($A\gtrace gq 0$ and $\alpha\gtrace gq 0$) if the $g$-trace
of the restriction of $A$ to any $q$-dimensional vector subspace
of $\cal V$ is nonnegative (respectively, positive); in other
words, for any choice of orthonormal vectors $e_1,\dots,e_q$ in
$\cal V$, we have $\sum A(e_j,e_j)\geq 0$ (respectively, $>0$).
We will apply the following elementary fact:
\begin{lem}\label{positive forms linear alg lemma}
Let $\left(\cal V,J\right)$ be a \cpx vector space of dimension
$n>1$, let~$g$ be a Hermitian inner product on $\cal V$ with
associated real skew-symmetric $(1,1)$-form~$\eta$, let
$\alpha$~and~$\beta$ be skew-symmetric real forms of type~$(1,1)$
on $\cal V$, and let \(\gamma\equiv\alpha\wedge\beta\wedge\eta^{n-2}\).
\begin{enumerate}

\item[(a)] If $\alpha\geq 0$ and $\beta\geqtrace g{n-1} 0$, then
$\gamma$ is nonnegative
(i.e. $\gamma/\eta^n\geq 0$).

\item[(b)] If $\alpha\geq 0$, $\beta\geq 0$, and $\gamma=0$, then
$\alpha\wedge\beta=0$.

\item[(c)] If $\alpha\geq 0$, $\beta\gtrace g{n-1}0$, and
$\gamma=0$, then $\alpha=0$.

\item[(d)] If $\alpha=\sqrt{-1}\omega_1\wedge\overline{\omega_1}$
for some form $\omega_1$ of type~$(1,0)$, $\beta\geq 0$, and
$\gamma=0$, then $\omega_1\wedge\beta=0$.

\item[(e)] If $\alpha=\sqrt{-1}\omega_1\wedge\overline{\omega_1}$
and $\beta=\sqrt{-1}\omega_2\wedge\overline{\omega_2}$ for some
pair of forms $\omega_1$ and $\omega_2$ of type~$(1,0)$ and
$\gamma=0$, then $\omega_1\wedge\omega_2=0$.

\end{enumerate}
\end{lem}
\begin{rmks}1. We used Part~(e) (which is easy to verify directly) in the proof
of Theorem~\ref{Cup product lemma general thm intro}.

\noindent 2. Parts~(b) and (d) may fail if we only assume that
$\beta\geqtrace g{n-1}0$ (in place of $\beta\geq 0$).
\end{rmks}
\begin{pf}
Corresponding to any \cpx basis $e_1,\dots,e_n$ for~$\cal V$, we
have
\begin{enumerate}

\item[(i)] The real basis $e_1,f_1=Je_1,\dots,e_n,f_n=Je_n$;

\item[(ii)] The real dual basis \(u_1,v_1=-u_1\circ
J,\dots,u_n,v_n=-u_n\circ J\); and

\item[(iii)] The \cpx basis
\(\zeta_1=u_1+\sqrt{-1}v_1,\dots,\zeta_n=u_n+\sqrt{-1}v_n,
\overline{\zeta_1},\dots,\overline{\zeta_n} \) for $\cal
V^*_{\C}$.

\end{enumerate}
We may choose the basis so that
\[
\eta=\sqrt{-1}\sum\zeta_j\wedge\overline{\zeta_j},\quad\alpha=\sqrt{-1}\sum
A_{ij}\zeta_i\wedge\overline{\zeta_j},\text{ and
}\beta=\sqrt{-1}\sum \lambda_j\zeta_j\wedge\overline{\zeta_j};
\]
where
\[
A_{ij}=\overline{A_{ji}}\quad\forall\,i,j\qquad\text{and}
\qquad\lambda_1\leq\lambda_2\leq\cdots\leq\lambda_n
\]
Thus
\begin{align*}
\gamma&=(\sqrt{-1})^{n-2}(n-2)!\sum_{1\leq i<j\leq
n}\alpha\wedge\beta\wedge\zeta_1\wedge\overline{\zeta_1}\wedge\cdots\wedge
\what{\zeta_i\wedge\overline{\zeta_i}}\wedge\cdots\wedge
\what{\zeta_j\wedge\overline{\zeta_j}}\wedge\cdots\wedge\zeta_n\wedge\overline{\zeta_n}\\
&=(\sqrt{-1})^n(n-2)!\sum_{1\leq i<j\leq
n}(A_{ii}\lambda_j+A_{jj}\lambda_i)
\zeta_1\wedge\overline{\zeta_1}\wedge\cdots\wedge\zeta_n\wedge\overline{\zeta_n}\\
&=2^n(n-2)!\sum_{i=1}^n\left[A_{ii}\left(\sum_{j\neq
i}\lambda_j\right)\right]u_1\wedge v_1\wedge\cdots\wedge u_n\wedge
v_n.
\end{align*}
Part~(a) follows immediately.

If $\alpha\geq 0$, $\beta\geq 0$, and $\gamma=0$, then $A_{ii}=0$
whenever $\lambda_k>0$ for some~$k\neq i$. Furthermore, whenever
$A_{ii}=0$, we have $A_{ij}=0$ for all~$j$ (for example, by the
Schwarz inequality).
Hence
\[
\alpha\wedge\beta=-\sum_{i,j=1}^n\sum_{k\neq
i,j}A_{ij}\lambda_k\zeta_i\wedge\overline{\zeta_j}\wedge\zeta_k\wedge\overline{\zeta_k}=0
\]
as claimed in~(b). Similar arguments give~(c).

Under the conditions in~(d), we have $A_{ij}=a_i\overline{a_j}$
where $\omega_1=\sum a_j\zeta_j$, and so $a_i=0$ whenever
$\lambda_k>0$ for some $k\neq i$. Hence
\[
\omega_1\wedge\beta=\sqrt{-1}\sum_{i=1}^n\sum_{k\neq
i}a_i\lambda_k\zeta_i\wedge\zeta_k\wedge\overline{\zeta_k}=0
\]
as claimed.

Under the conditions in~(e), we get $A_{ij}=a_i\overline{a_j}$ for
all~$i,j$ and $0=\lambda_1=\dots=\lambda_{n-1}\leq\lambda_n=|b|^2$
where $\omega_1=\sum a_j\zeta_j$ and $\omega_2=b\zeta_n$. Hence
$a_1=\dots=a_{n-1}=0$ if $b\neq 0$ and, therefore,
\[
\omega_1\wedge\omega_2=\sum_{i=1}^{n-1}a_ib\zeta_i\wedge\zeta_n=0
\]
as claimed.
\end{pf}

It will also be convenient to fix a \cinf \fn $\chi\colon\R\to\R$
\st $\chi'\geq 0$ and $\chi''\geq 0$ on~$\R$, $\chi(t)=0$ for
$t\leq 0$, and $\chi(t)=t-1$ for $t\geq 2$.

\begin{thm}\label{cup product lemma very general thm}
Let $\vphi$ be a real-valued \cinf \fn on $X$, let
$\alpha=2\sqrt{-1}\partial\dbar\vphi=dd^c\vphi$, and let $\beta$
be a \cinf real form of type~$(1,1)$ on~$X$. Assume that
\begin{enumerate}
\item[(i)] The real $(n,n)$-form
$\gamma\equiv\alpha\wedge\beta\wedge\eta^{n-2}$ satisfies
$\gamma\geq 0$;

\item[(ii)] We have $d^c\vphi\wedge
d\left(\beta\wedge\eta^{n-2}\right)\equiv 0$; and

\item[(iii)] For some point $p\in X$, we have
\[
\liminf_{r\to\infty}\frac 1r\int_{B(p;2r)\sm
B(p;r)}|d\vphi|_g|\beta|_g\,dV_g=0.
\]

\end{enumerate}
Then the following hold:
\begin{enumerate}

\item[(a)] We have $\gamma\equiv 0$ on $X$.

\item[(b)] At any point~$x\in X$ at which both $\alpha_x\geq 0$
and $\beta_x\geq 0$ hold, we have \((\alpha\wedge\beta)_x=0\).

\item[(c)] At any point~$x\in X$ at which both $\alpha_x\geq 0$
and $\beta_x=\sqrt{-1}\omega\wedge\bar\omega$ for some
$\omega\in\left(T_x^{1,0}X\right)^*$, we have
$\alpha_x\wedge\omega=0$.

\item[(d)] If
$\sqrt{-1}\partial\vphi\wedge\dbar\vphi\wedge\beta\wedge\eta^{n-2}\geq
0$, then
$\partial\vphi\wedge\dbar\vphi\wedge\beta\wedge\eta^{n-2}\equiv
0$.

\item[(e)] If
$\sqrt{-1}\partial\vphi\wedge\dbar\vphi\wedge\beta\wedge\eta^{n-2}\geq
0$, then, at any point~$x\in X$ at which $\beta_x\geq 0$ holds, we
have $(\partial\vphi\wedge\beta)_x=0$.

\item[(f)] If
$\sqrt{-1}\partial\vphi\wedge\dbar\vphi\wedge\beta\wedge\eta^{n-2}\geq
0$, then, at any point $x\in X$ at which
$\beta_x=\sqrt{-1}\omega\wedge\bar\omega$ for some
$\omega\in\left(T_x^{1,0}X\right)^*$, we have
$(\partial\vphi)_x\wedge\omega=0$.

\end{enumerate}
\end{thm}
\begin{rmks}
1. By Lemma~\ref{positive forms linear alg lemma}, the
condition~(i) holds if, for example, at each point
$x\in\supp\alpha\cap\supp\beta$ we have $\alpha_x\geq 0$ and
$\beta_x\geqtrace g{n-1}0$ or we have $\alpha_x\geqtrace g{n-1}0$
and $\beta_x\geq 0$.

\noindent 2. Condition~(ii) holds if, for example, $g$ is K\"ahler
and $\beta$ is closed off of the set of critical points
of~$\vphi$.

\noindent 3. The condition
$\sqrt{-1}\partial\vphi\wedge\dbar\vphi\wedge\beta\wedge\eta^{n-2}\geq
0$ (in Parts~(d)--(f)) holds if, for example, at each point
$x\in\supp(d\vphi)\cap\supp\beta$, we have $\beta_x\geqtrace
g{n-1}0$.

\noindent 4. The condition~(iii) holds if, for example, $\beta\in
L^1$ and, for some constant $C>0$, $|d\vphi|_g\leq C(r+1)$ on
$B(p;r)$ for every $r>0$.

\noindent 5. Clearly, Theorem~\ref{Cup product lemma general thm
intro} is a special case of the above theorem (one simply takes $\vphi=\rho$ where
$\rho$ is a \plh \fn with $\real(\omega_1)=d\rho$).
\end{rmks}
\begin{pf*}{Proof of Theorem~\ref{cup product lemma very general thm}} We have
\[
\gamma=dd^c\vphi\wedge\beta\wedge\eta^{n-2}=d\left[d^c\vphi\wedge\beta\wedge\eta^{n-2}\right].
\]
Hence, for each $r>0$, Stokes' theorem gives
\[
\int_X\tau_r\gamma=-\int_{B(p;2r)\sm B(p;r)}d\tau_r\wedge
d^c\vphi\wedge\beta\wedge\eta^{n-2}.
\]
On the other hand, for some constant $C=C(n)>0$, we have
\[
\left|d\tau_r\wedge
d^c\vphi\wedge\beta\wedge\eta^{n-2}\right|_g\leq \frac
Cr|d\vphi|_g|\beta|_g.
\]
Thus, for a suitable sequence $\seq r\nu$ with $r_\nu\to\infty$,
we get
\[
\int_X\gamma\leftarrow\int_X\tau_{r_\nu}\gamma\leq\frac
C{r_\nu}\int_{B(p;2r_\nu)\sm B(p;r_\nu)}|d\vphi|_g|\beta|_gdV_g\to
0.
\]
Since $\gamma\geq 0$, we get $\gamma\equiv 0$ as claimed in~(a)
and Lemma~\ref{positive forms linear alg lemma} gives~(b)~and~(c)
as well.

Assume now that
$\sqrt{-1}\partial\vphi\wedge\dbar\vphi\wedge\beta\wedge\eta^{n-2}\geq
0$ on~$X$. Given a constant $a\in\R$, the nonnegative \cinf
\fnns~$\psi=\chi(\vphi-a)$ satisfies
\[
dd^c\psi=2\sqrt{-1}\partial\dbar\psi=\chi'(\vphi-a)\alpha
+2\sqrt{-1}\chi''(\vphi-a)\partial\vphi\wedge\dbar\vphi
\]
and
\begin{align*}
dd^c\psi^2&=d(2\psi d^c\psi)=2\sqrt{-1}\partial\dbar\psi^2\\
&=4\psi\sqrt{-1}\partial\dbar\psi+4\sqrt{-1}\partial\psi\wedge\dbar\psi\\
&=2\chi(\vphi-a)\chi'(\vphi-a)\alpha
+\left(\chi(\vphi-a)\chi''(\vphi-a)+\left[\chi'(\vphi-a)\right]^2\right)4\sqrt{-1}\partial\vphi\wedge\dbar\vphi
\end{align*}
Thus, by the hypotheses, \(d\theta=\xi\equiv
dd^c\psi^2\wedge\beta\wedge\eta^{n-2}\geq 0\), where
\(\theta\equiv 2\psi d^c\psi\wedge\beta\wedge\eta^{n-2}\). For
every regular value $R>1$ of $\psi$, let $\xi_R\geq 0$ be the
product of $\xi$ and the characteristic \fn of the set
\(\setof{x\in X}{\psi\leq R}\), let $\psi_R=\cal
M_R(\psi)=\min(\psi,R)$, and let $\theta_R$ be the Lipschitz \cont
form given by
\[
\theta_R\equiv2\psi_Rd^c\psi\wedge\beta\wedge\eta^{n-2}.
\]
On the set $\setof{x\in X}{\psi(x)>1}\supset\setof{x\in
X}{\psi(x)\geq R}$, we have $\psi=\vphi-1$ and hence
\[
dd^c\psi\wedge\beta\wedge\eta^{n-2}=\gamma=0.
\]
Therefore $\xi_R=d\theta_R$ on $X\sm\psi\inv(R)$. For each $r>0$,
Stokes' theorem gives
\[
\int_X\tau_r\xi_R=-\int_{B(p;2r)\sm B(p;r)}d\tau_r\wedge\theta_R.
\]
On the other hand, since $0\leq\chi'\leq 1$, we have
\(|d\tau_r\wedge\theta_R|_g\leq 2Cr\inv R|d\vphi|_g|\beta|_g\)
(where $C=C(n)>0$ as before). Letting $r=r_\nu\to\infty$ for a
suitable sequence $\seq r\nu$, we get
\[
\int_X\xi_R=0.
\]
Since $\xi_R\geq 0$, we must have $\xi_R\equiv 0$ and, letting
$R\to\infty$, we get (since $\gamma\equiv 0$)
\begin{align*}
0=\xi&=dd^c\psi^2\wedge\beta\wedge\eta^{n-2}\\
&=\left(\chi(\vphi-a)\chi''(\vphi-a)+\left[\chi'(\vphi-a)\right]^2\right)
4\sqrt{-1}\partial\vphi\wedge\dbar\vphi\wedge\beta\wedge\eta^{n-2}.
\end{align*}
Now, given a point $x\in X$, we may choose $a\in\R$ so that
$2<\vphi(x)-a$. Hence
\[
0=-4\inv\sqrt{-1}\xi_x=\left(\partial\vphi\wedge\dbar\vphi\wedge\beta\wedge\eta^{n-2}\right)_x,
\]
as claimed in~(d), and Lemma~\ref{positive forms linear alg lemma}
gives~(e)~and~(f).
\end{pf*}

We have the following two immediate consequences:
\begin{cor}\label{cup prod lemma plsh fn and form lem}
Let $\vphi$ be a \cinf \plsh \fnns, let $Z$ be the set of critical
points of~$\vphi$, and let $\beta$ be a \cinf real form of
type~$(1,1)$ on~$X$. Assume that $\beta\restrict{X\sm Z}$ is
closed and nonnegative, $g\restrict{X\sm Z}$ is K\"ahler, and, for
some point $p\in X$,
\[
\liminf_{r\to\infty}\frac 1r\int_{B(p;2r)\sm
B(p;r)}|d\vphi|_g|\beta|_g\,dV_g=0.
\]
Then
\[
\partial\dbar\vphi\wedge\beta\equiv
0\qquad\text{and}\qquad\partial\vphi\wedge\beta\equiv 0.
\]
\end{cor}

\begin{cor}\label{cup prod lemma plsh fn and holo 1form lem}
Let $\vphi$ be a \cinf \plsh \fn on~$X$, let $Z$ be the set of
critical points of~$\vphi$, and let $\omega$ be a \cinf form of
type $(1,0)$ on~$X$ \st $\omega\restrict{X\sm Z}$ is closed (hence
\holons), $g\restrict{X\sm Z}$ is K\"ahler, and, for some point
$p\in X$, we have
\[
\liminf_{r\to\infty}\frac 1r\int_{B(p;2r)\sm
B(p;r)}|d\vphi|_g|\omega|^2_g\,dV_g=0.
\]
Then
\[
\partial\dbar\vphi\wedge\omega\equiv
0\qquad\text{and}\qquad\partial\vphi\wedge\omega\equiv 0.
\]
\end{cor}

\begin{rmk}
The above limit inferior is~$0$ if, for example, $\omega$ is
in~$L^2$ and, for some~$C>0$,  $|d\vphi|_g\leq C(r+1)$ on
$B(p;2r)$ for each $r>0$.
\end{rmk}

\begin{defn}
Let $q$ be a positive integer. A \cinf real-valued \fnns~$\vphi$
on an open subset~$\Omega$ of~$X$ is of class
$\plshclass^\infty(g,q)$ (of class $\strplshclass^\infty(g,q)$) if
$\sqrt{-1}\partial\dbar\vphi\geqtrace gq0$ (respectively,
$\sqrt{-1}\partial\dbar\vphi\gtrace gq 0$).
\end{defn}
This class of \fns was first introduced by Grauert and
Riemenschneider~\cite{Grauert-Riemenschneider} and has since been
applied in several contexts (see, for example, \cite{Siu-complex
analyticity of harmonic maps}, \cite{Wu-On certain Kahler
manifolds which are q-complete}, \cite{NR-BH Regular hyperbolic
Kahler}, \cite{Joita Traces of cvx domains},
\cite{Fraboni-Covering cpx mfld}). Theorem~\ref{cup product lemma
very general thm} immediately gives the following:
\begin{cor}\label{cup product for n-1plsh and holo 1form cor}
Let $\vphi\in\plshclass^\infty(g,n-1)(X)$, let
$Z$ be the set of critical points of~$\vphi$, and let $\omega$ be
a \cinf form of type~$(1,0)$ on $X$ \st $\omega\restrict{X\sm Z}$
is closed (hence \holons), $g\restrict{X\sm Z}$ is K\"ahler, and,
for some point $p\in X$, we have
\[
\liminf_{r\to\infty}\frac 1r\int_{B(p;2r)\sm
B(p;r)}|d\vphi|_g|\omega|^2_g\,dV_g=0.
\]
Then $\partial\vphi\wedge\omega\equiv 0$ on $X$.
\end{cor}

Theorem~\ref{cup product plh fn and form on superlevel thm} may be
considered as a consequence of Corollary~\ref{cup prod lemma plsh
fn and holo 1form lem} (or Corollary~\ref{cup product for n-1plsh
and holo 1form cor}) and Theorem~\ref{L2 Castel de Franch on open
set thm } by fixing a number $b$ with $a<b<\sup_\Omega\rho_1$ and
setting
\[
\vphi=\left\{
\begin{aligned}
\chi\left(\frac{2\rho_1-2a}{b-a}\right)&\qquad\text{on }\Omega\\
0&\qquad\text{on }X\sm\Omega.
\end{aligned}
\right.
\]

Similar arguments also give the following theorem which may be
viewed both as a variant of Lemma~2.6 of \cite{NR Filtered
ends} (see also Lemma~2.1 of \cite{NR-BH Regular hyperbolic
Kahler}) and of Lemma~2.7 of \cite{NR Filtered ends} and which,
in the bounded geometry complete K\"ahler case, is a
generalization of both.

\begin{thm}\label{cup product for band thm}
Let $\rho_1$ and $\rho_2$ be two real-valued \plh \fns on a
nonempty domain~$Y$ in $X$. Assume that, for some pair of
constants $a,b$ with $\inf \rho _1 <a<b<\sup\rho _1$ and some
\comp $\Omega$ of $\setof{x\in Y}{a<\rho _1(x)<b}$, we have the
following:
\begin{enumerate}
\item[(i)] $\overline\Omega\subset Y$;

\item[(ii)] The metric~$g\restrict\Omega$ is K\"ahler;

\item[(iii)] The form $d\rho_2\restrict\Omega$ is bounded; and

\item[(iv)] $\int_\Omega|d\rho_j|^2_g\,dV_g <\infty$ for $j=1,2$.

\end{enumerate}
Then $\partial\rho_1\wedge\partial\rho_2\equiv 0$ on $Y$.
Furthermore, if $d\rho_1$ and $d\rho_2$ are linearly independent
and $(X,g)$ has bounded geometry along $\Omega$, then there exist
a surjective proper \holo mapping $\Phi\colon\Omega\to S$ of
$\Omega$ onto a Riemann surface~$S$ with $\Phi_*\ol_\Omega=\ol_S$
and \plh \fns $\alpha_1$ and $\alpha_2$ on $S$ \st
$\rho_j\restrict\Omega=\Phi^*\alpha_j$ for $j=1,2$.
\end{thm}
\begin{pf}
The coarea formula gives
\[
\int_a^b\left[\int_{\rho_1\inv(t)\cap\Omega}|d\rho_1|_g\,d\sigma_g\right]\,dt
=\int_\Omega|d\rho_1|_g^2\,dV_g<\infty.
\]
Thus we may choose regular values $A$ and $B$ of $\rho_1$
in~$\rho_1(\Omega)$ \st $a<A<B<b$ and \stns, for
$M=\rho_1\inv(A)\cap\Omega$ and $N=\rho_1\inv(B)\cap\Omega$, we
have
\[
\int_M|d\rho_1|_g\,d\sigma_g<\infty\qquad\text{and}\qquad
\int_N|d\rho_1|_g\,d\sigma_g<\infty.
\]
In particular, we have
$\Theta=\setof{x\in\Omega}{A<\rho_1(x)<B}\neq\emptyset$ and
$\overline\Theta\subset\Omega\subset Y$.
The real $(n,n)$-form
\[
\gamma\equiv(\sqrt{-1}\partial\rho_1\wedge\dbar\rho_1)
\wedge(\sqrt{-1}\partial\rho_2\wedge\dbar\rho_2)\wedge\eta^{n-2}
\]
on $\Omega$ then satisfies
\[
0\leq\gamma=\frac 1{16}dd^c(\rho_1)^2\wedge
dd^c(\rho_2)^2\wedge\eta^{n-2}=\frac 14d(\rho_1d^c\rho_1)\wedge
d(\rho_2d^c\rho_2)\wedge\eta^{n-2}=d\theta;
\]
where $\theta$ is the $L^1$ (and $L^2$) form on~$\Omega$ given by
\[
\theta=\frac 14\rho_1d^c\rho_1\wedge
d(\rho_2d^c\rho_2)\wedge\eta^{n-2}=\frac 12\rho_1d^c\rho_1\wedge(\sqrt{-1}
\partial\rho_2\wedge\dbar\rho_2)\wedge\eta^{n-2}.
\]

We may form a complete Hermitian metric $h$ in $\Omega$ \st $h\geq
g$ on $\Omega$ and $h=g$ on a \nbd of $\overline\Theta$ and,
fixing a point $p\in\Theta$ and applying the Gaffney
construction~\cite{Gaffney} and \cinf approximation, we get a
collection of nonnegative \cinf \fns $\seq\kappa r_{r>0}$ \stns,
for each $r>0$, we have $0\leq\kappa_r\leq 1$ on $X$,
$\supp\kappa_r\subset B_h(p;2r)\Subset\Omega$, $\kappa_r\equiv 1$
on $B_h(p;r)$, and $|d\kappa_r|_h\leq 2/r$.

For each $r>0$, Stokes' theorem gives
\begin{align*}
\int_\Theta\kappa_r\gamma&=\int_N\kappa_r\theta-\int_M\kappa_r\theta-\int_{[B_h(p;2r)\sm
B_h(p;r)]\cap\Theta}d\kappa_r\wedge\theta\\
&=\frac B4\int_N\kappa_rd^c\rho_1\wedge
d(\rho_2d^c\rho_2)\wedge\eta^{n-2}-\frac
A4\int_M\kappa_rd^c\rho_1\wedge d(\rho_2d^c\rho_2)\wedge\eta^{n-2}\\
&\qquad\qquad-\int_{[B_h(p;2r)\sm
B_h(p;r)]\cap\Theta}d\kappa_r\wedge\theta\\
&=\frac B4\int_Nd\kappa_r\wedge d^c\rho_1\wedge
\rho_2d^c\rho_2\wedge\eta^{n-2}-\frac
A4\int_Md\kappa_r\wedge d^c\rho_1\wedge\rho_2d^c\rho_2\wedge\eta^{n-2}\\
&\qquad\qquad-\int_{[B_h(p;2r)\sm
B_h(p;r)]\cap\Theta}d\kappa_r\wedge\theta
\end{align*}
Since $|d\rho_2|_h\leq|d\rho_2|_g$, $|d\rho_2|_h$ is bounded on~$\Omega$ and hence
$\rho_2$ must have at most linear growth; i.e.~for some constant
$C>0$, we have $|\rho_2|\leq C(r+1)$ on $B_h(p;r)$ for
every~$r>0$. Moreover, $|d^c\rho_1|_g$ is in~$L^1$ on $M$ and $N$
by the choice of $A$ and $B$. Thus we may apply the dominated
convergence theorem as $r\to\infty$ to get
\[
\int_\Theta\gamma=0.
\]
Since $\gamma\geq 0$, we get $\gamma\equiv 0$ on~$\Theta$ and
Lemma~\ref{positive forms linear alg lemma} implies that
$\partial\rho_1\wedge\partial\rho_2\equiv 0$ on $\Theta$ and,
therefore, on~$Y$.

If $d\rho_1$ and $d\rho_2$ are linearly independent and $(X,g)$
has bounded geometry along $\Omega$, then we may apply
Corollary~\ref{L2 Castelnuovo on open set for partial exact forms
cor} to get the desired proper \holo mapping to a Riemann surface
$\Phi\colon\Omega\to S$.
\end{pf}

In particular, we get the following weak version of Lemma~2.1
of~\cite{NR-BH Regular hyperbolic Kahler} (see also Lemma~2.6
of~\cite{NR Filtered ends}) which suffices for the proof of the
main result of~\cite{NR Filtered ends}:
\begin{cor}\label{cup product cpt fiber plh cor}
Let $\rho_1$ and $\rho_2$ be two real-valued \plh \fns on $X$.
Assume that $\rho_1$ has a nonempty \cpt fiber $F$ for which $g$
is K\"ahler on a \nbd of~$F$. Then
$\partial\rho_1\wedge\partial\rho_2\equiv 0$ on $X$. Furthermore,
if $d\rho_1$ and $d\rho_2$ are linearly independent, then there
exist a proper \holo mapping $\Phi\colon\Omega\to S$ of a
\nbdns~$\Omega$ of $F$ in~$X$ onto a Riemann surface~$S$ with
$\Phi_*\ol_\Omega=\ol_S$ and \plh \fns $\alpha_1$ and $\alpha_2$
on $S$ \st $\rho_j\restrict\Omega=\Phi^*\alpha_j$ for $j=1,2$.
\end{cor}

\section{$L^2$ Castelnuovo-de Franchis for an end}\label{L2 Castelnuovo for an end sect}

In this section, the arguments of Section~\ref{L2 Castel de
Franchis sect} are extended in order to obtain the following
version of the $L^2$~Castelnuovo-de Franchis theorem which is
applied in~\cite{NR-BH bound geom hyperbolic Kahler} and is
also of separate interest:

\begin{thm}\label{L2 Castelnuovo for end thm}
Let $(X,g)$ be a \con non\cpt complete Hermitian \mfldns, let $E$
be a special end of type~(BG) in~$X$, and let $\omega_1$ and
$\omega_2$ be linearly independent closed \holo $1$-forms on~$E$
\st $g\restrict E$ is K\"ahler, $\omega_1$ is in~$L^2$, and
$\omega_1\wedge\omega_2\equiv 0$ on~$E$. Then there exist a
surjective proper \holo mapping $\Phi\colon\Omega\to S$ of a
nonempty open subset $\Omega$ of~$E$ onto a Riemann surface~$S$
with $\Phi_*\ol_\Omega=\ol_S$ and \holo $1$-forms $\theta_1$ and
$\theta_2$ on $S$ \st $\omega_j\restrict\Omega=\Phi^*\theta_j$ for $j=1,2$.
\end{thm}
\begin{rmk}
In particular, by Lemma~\ref{Open set to Riemann surface gives
global lemma}, if, in addition, $g$ is K\"ahler and $X$ admits a
special ends decomposition, then $X$ admits a proper \holo mapping
onto a Riemann surface.
\end{rmk}

\begin{cor}\label{L2 Castelnuovo for end plh fn cor}
Let $(X,g)$ be a \con non\cpt complete Hermitian \mfldns, let $E$
be a special end of type~(BG) in~$X$, and let
$\rho_1$~and~$\rho_2$ be two real-valued \plh \fns on~$E$ \st
$g\restrict E$ is K\"ahler, $d\rho_1$ and $d\rho_2$ are linearly
independent, $\rho_1$ has finite energy, and
$\partial\rho_1\wedge\partial\rho_2\equiv 0$ on~$E$. Then there
exist a surjective proper \holo mapping $\Phi\colon\Omega\to S$ of
a nonempty open subset $\Omega$ of~$E$ onto a Riemann surface~$S$
with $\Phi_*\ol_\Omega=\ol_S$ and real-valued \plh \fns $\alpha_1$
and $\alpha_2$ on $S$ \st $\rho_j\restrict\Omega=\Phi^*\alpha_j$ for $j=1,2$.

In particular, if there exists a nonconstant \holo \fn with finite
energy on $E$, then there exists a surjective proper \holo mapping
$\Phi\colon\Omega\to S$ of a nonempty open subset $\Omega$ of~$E$
onto a Riemann surface~$S$ with $\Phi_*\ol_\Omega=\ol_S$.
\end{cor}

For the proof of Theorem~\ref{L2 Castelnuovo for end thm}, we will
show that, under the assumption that the associated \holo mapping
to $\mathbb P^1$ has no nonempty \cpt levels, one gets a \holo \fn
which is defined in $E$ outside a \rel \cpt \nbd of $\partial E$
and which vanishes at infinity.  In particular, the end is then
parabolic by the following observation of J.~Wang (cf.~Lemma~1.3
of~\cite{NR-BH Regular hyperbolic Kahler}):
\begin{lem}\label{holo fn vanish at infinity parabolic lemma}
Let $(X,g)$ be a complete Hermitian manifold and let $E$ be an end
in~$X$ \st $g\restrict E$ is K\"ahler. Assume that there exists a
nonconstant \holo \fnns~$h$ on a \nbd of~$\overline E$ \st
$\lim_{x\to\infty}h\restrict{\overline E}(x)=0$. Then $E$ is a
parabolic end.
\end{lem}
\begin{pf}
The \fn $\vphi =-\log |h|^2\colon\overline E\to(-\infty,\infty]$
is superharmonic on~$E$ and $\vphi(x)\to\infty$ as $x\to\infty$ in
$\overline E$. In particular, we may assume that $\vphi$ is
positive. Suppose $\alpha$ is a nonnegative bounded sub\harm \fn
on~$X$ which vanishes on $X\sm E$. Given $\epsilon>0$, we have
$\epsilon\vphi>0=\alpha$ on $\partial E$ and
$\epsilon\vphi>\sup\alpha$ on the complement in~$E$ of a
sufficiently large \cpt subset of $X$. It follows that
$0\leq\alpha<\epsilon\vphi$ on~$\overline E$ for every
$\epsilon>0$ and, therefore, that $\alpha\equiv 0$. Thus $E$ does
not admit an admissible sub\harm \fn and hence $E$ is a parabolic
end.
\end{pf}

We will apply the following lemma which is a consequence of the
work of Grauert and Riemenschneider~\cite{Grauert-Riemenschneider}
(for a \rel \cpt domain~$E$), of Gromov~\cite{Gro-Kahler
hyperbolicity} and of Li~\cite{Li Structure complete Kahler} (for
$E=X$), and of Siu~\cite{Siu-complex analyticity of harmonic
maps} (for a harmonic mapping of a \rel \cpt domain into a manifold
satisfying certain curvature conditions).
\begin{lem}[Grauert-Riemenschneider, Li, Siu (see Lemma~3.2 of~\cite{NR-BH Regular hyperbolic
Kahler})]\label{Grauert-Riemens-Li-Siu on n-1 plsh boundary lemma}
Let $(X,g)$ be a \con complete Hermitian manifold of dimension
$n$, let $E$ be a (not necessarily \rel \cptns) domain with smooth
\cpt (possibly empty) boundary in~$X$, let $\vphi$ be \cinf
real-valued \fn on~$X$ \st $d\vphi\neq 0$ at every point in
$\partial E$ and \st $E=\setof{x\in X}{\vphi (x) <0}$, and, for
each point $x\in\partial E$, let
\[
\tau (x)=\text{{\rm
tr}}\,\left[\lev\vphi\restrict{T_x^{1,0}(\partial E)}\right].
\]
Assume that $g\restrict E$ is K\"ahler and that $\tau \geq 0$ on
$\partial E$. Then we have the following:
\begin{enumerate}
\item[(a)] If $\beta $ is a~$\cinf $ \fn on~$\overline E$ \st
$\beta $ is \harm on $E$, $\beta $ satisfies the tangential
Cauchy-Riemann equation $\dbar _b\beta =0$ on $\partial E$, and
there is a sequence of positive real numbers $R_m\to\infty$ and a
point~$p\in X$ \st
\[
\lim_{m\to\infty}\frac 1{R_m^2}\|\nabla\beta\|
^2_{L^2(B_p(R_m)\cap E)}=\lim_{m\to\infty}\frac
1{R_m^2}\int_{B_p(R_m)\cap E}|\nabla\beta |^2\,dV=0,
\]
then $\beta$ is \plh on $E$.

\item[(b)] If $E$ is a hyperbolic end of~$X$, then $\tau\equiv 0$
on $\partial E$.
\end{enumerate}
\end{lem}

As suggested by the above lemma, functions of class
$\strplshclass^\infty(g,q)$ will play a role in the proof of the
parabolic case of the theorem. The following fact is contained implicitly in the
work of Richberg~\cite{Richberg}, Greene and Wu~\cite{Green-Wu
Embedding of open}, Ohsawa~\cite{Ohsawa completeness of noncpt
analytic}, Col\c toiu~\cite{Coltoiu pluripolar}, and
Demailly~\cite{Demailly Cohomology of q-cvx} (see~\cite{NR-BH
Weakly 1-complete}):
\begin{prop}[Richberg, Greene-Wu, Ohsawa,
Col\c toiu, Demailly]\label{qplsh fn on nbd prop} Suppose $(X,g)$
is a Hermitian manifold of dimension~$n>1$ and $Y$ is a nowhere
dense \anal subset with no \cpt \ircompsns. Then there exists a
\cinf \exh \fn $\vphi$ on $X$ which is of class
$\strplshclass^\infty(g,n-1)$ on a \nbd of $Y$ in~$X$.
\end{prop}
The following proposition will give the parabolic case of the
theorem and is also of separate interest:

\begin{prop}\label{parabolic end vanish at infinity end prop}
Let $(X,g)$ be a \con non\cpt complete Hermitian \mfld and let $E$
be a special end of type~(BG) in~$X$. Assume that $g\restrict E$
is K\"ahler and that there exists a nonconstant \holo \fnns~$h$
on~$E$ which vanishes at infinity in~$X$. Then there exists a
proper \holo mapping of some nonempty open subset of $E$ onto a
Riemann surface.
\end{prop}
\begin{rmks}
1. The end $E$ is necessarily parabolic by Lemma~\ref{holo fn
vanish at infinity parabolic lemma}.

\noindent 2. The above proposition also holds if $E$ is special of
type~(W) instead of type~(BG); as the proof below together with
Proposition~2.3 of~\cite{NR-BH Weakly 1-complete} shows.
\end{rmks}
\begin{pf*}{Proof of Proposition~\ref{parabolic end vanish at infinity end prop}} We may assume without loss of generality that $n=\dim
X>1$. The idea of the proof is to construct an end with a defining
\fn of class~$\plshclass^\infty(g,n-1)$ at the \cpt boundary
and to apply Lemma~\ref{Grauert-Riemens-Li-Siu on n-1 plsh
boundary lemma}. Let $Z=h\inv(0)\subset E$ and fix open subsets
$\Omega_1$,~$\Omega_2$,~and~$\Omega_3$ \st
\[
\partial E\subset\Omega_1\Subset\Omega_2\Subset\Omega_3\Subset X
\]
and $Z\cap\Omega_3\sm\overline\Omega_1$ is empty or has no \cpt
\ircomps (one can form these sets by choosing open sets
$\Omega_1\Subset\Omega_3'\Subset X$ and setting $\Omega_3=\Omega_3'\sm F$,
where $F$ is a finite subset of $E\sm\overline{\Omega_1}$ which
contains a point in $Z_0\cap\Omega_3'\sm\overline{\Omega_1}$ for
each \ircompns~$Z_0$ of $Z$ which meets
$\Omega'_3\sm\overline{\Omega_1}$). By Proposition~\ref{qplsh fn
on nbd prop} (Richberg, Greene-Wu, Ohsawa, Col\c toiu, Demailly),
there exists a \nbdns~$V$ of $Z\cap\Omega_3\sm\overline{\Omega_1}$
in $E\cap\Omega_3\sm{\overline\Omega_1}$ and a positive \fn
$\psi\in\strplshclass^\infty(g,n-1)(V)$ which exhausts
$Z\cap\Omega_3\sm\overline{\Omega_1}$. We may also fix constants
$a$~and~$b$ with
\[
a>b>\max_{Z\cap\partial\Omega_2}\psi\quad(a>b>0\text{ if
}Z\cap\partial\Omega_2=\emptyset).
\]
After shrinking~$V$, we may assume that there is a \nbd of
$\partial(\Omega _3\sm\overline\Omega_1)$ in~$X$ \st $\psi>a$ at
points in~$V$ which lie in this \nbdns.
Therefore, for $\epsilon>0$ sufficiently small, the set
\[
\setof{x\in V}{|h(x)|\leq\epsilon\text{ and }\psi(x)\leq a}
\]
is \cpt (and possibly empty) and the \cpt (and possibly empty) set
\[
\setof{x\in\partial\Omega_2}{|h(x)|\leq\epsilon}
\]
is contained in $\setof{x\in V}{\psi(x)<b}$. Choosing a \cinf
nondecreasing convex \fn $\chi:\R\to\R$ which vanishes on the
interval $(-\infty,-\log(a-b)]$ and which approaches~$+\infty $ at
$+\infty$, we obtain a \fnns~$\vphi$ of
class~$\plshclass^\infty(g,n-1)$ on the open set
\[
\Omega=\setof{x\in E\sm\Omega_2}{|h(x)|<\epsilon}\cup\setof{x\in
V\cap\Omega_2}{|h(x)|<\epsilon\text{ and }\psi(x)<a}
\]
by defining
\[
\vphi(x)=\left\{
\begin{alignedat}{2}
&-\log(\epsilon^2-|h(x)|^2)&\quad\text{ if }x\in E\cap\Omega\sm\Omega _2\\
&\chi(-\log(a-\psi(x))-\log(\epsilon^2-|h(x)|^2)&\quad\text{ if
}x\in E\cap\Omega\cap\Omega _2
\end{alignedat}
\right.
\]
Moreover, $\vphi\geq-\log\epsilon^2$ on~$\Omega$, and, since $h$
vanishes at infinity in~$X$, $\Omega$ is not \rel \cpt in~$X$, the
boundary $\partial\Omega$ in~$X$ (and in~$E$) is a nonempty \cpt
subset of~$E$, $\vphi\to\infty$ at~$\partial \Omega$, and
$\vphi\to-\log\epsilon^2$ at infinity in~$X$. Therefore, if $c$ is
a regular value of~$\vphi$ with $c>-\log\epsilon^2$ and $E_0$ is a
\concomp of the set $\setof{x\in\Omega}{\vphi(x)<c}$ with non\cpt
closure in~$X$, then $E_0$ is a special end of type~(BG) in~$X$
with smooth nonempty (\cptns) boundary and $E_0$ admits a defining
\fn of class $\plshclass^\infty(g,n-1)$ (on a \nbd of
$\overline E_0$) with nonvanishing differential at each boundary
point.

As in the proof of Theorem~2.6 of~\cite{NR-Structure theorems},
one may apply a theorem of Nakai~\cite{Nakai1}, \cite{Nakai2} and
a theorem of Sullivan~\cite{Sul} to obtain a \cont \fn $\rho$
on~$\overline E_0$ \st $\rho$ is \harm on~$E_0$, $\rho=0$ on
$\partial E_0$, $\rho(x)\to\infty$ at infinity, and the $L^2$ norm
of $\grad\rho$ on a ball of radius~$R$ is equal to $o(R)$.
Therefore $\rho$ is \plh by Lemma~\ref{Grauert-Riemens-Li-Siu on
n-1 plsh boundary lemma} and $\rho$ has a (nonempty) \cpt fiber
in~$E_0$. The pair of \plh \fns $\rho$ and $\real h\restrict{E_0}$
must have linearly independent differentials on~$E_0$ (since
$\rho\to\infty$ while $\real h\to 0$ at infinity), so
Corollary~\ref{cup product cpt fiber plh cor} now gives the claim.
\end{pf*}

\begin{lem}\label{finite volume fibers in end lem}
Let $(X,g)$ be a \con non\cpt complete Hermitian \mfldns, let $E$
be an end of~$X$, let $h\colon E\to\mathbb P^1=\C\cup\set\infty$
be a \holo mapping with no nonempty \cpt levels, let
$C=\setof{x\in E}{(h_*)_x=0}$ be the set of critical points, let
$Z$ be the union of all \concomps~$A$ of~$C$ for which $\bar A$ is
non\cpt and meets~$\partial E$, and let $r\in(0,\infty)$. Assume
that $g\restrict E$ is K\"ahler, the extended real line $\mathbb
S^1=\R\cup\set\infty\subset\mathbb P^1$ is contained in $\mathbb
P^1\sm h(Z)$, $(\R\sm\set 0)\cup\set\infty\subset\mathbb P^1\sm
h(C)$ (the set of regular values), and
\[
h\inv(t)\Subset X\quad\forall\,t\in(0,r).
\]
Then, for some \rel \cpt open subset $\Omega$ of~$X$
containing~$\partial E$, the \fn
$$
t\mapsto\Vol\left(h\inv(t)\sm\overline\Omega\right)
$$
is bounded on the interval $(0,r)$.
\end{lem}
\begin{pf} Let $n=\dim X$ and let $\eta$ be the real $(1,1)$-form associated to~$g$.
The idea is to apply the Stokes theorem arguments of Stoll~[St]
while keeping track of the boundary integrals. The set $M=h\inv
(\R\cup\set\infty)$ is a (properly embedded) real \anal subset
of~$E$ and the subset $M\sm C$ is a (properly embedded) oriented
(with nonvanishing $(2n-1)$-form $(h^*d\theta )\wedge\eta^{n-1}$)
real \anal submanifold of dimension~$2n-1$ in~$E\sm C$ (assuming,
as we may, that $M$ is nonempty). We also have $M\cap C=M\cap
h\inv(0)\cap C\sm Z$.

We may choose a real-valued \cinf \fnns~$\rho$ on~$X$, a \rel \cpt
\nbd $W$ of $\partial E$ in~$X$, and a constant~$\epsilon>0$ \st
$[-\epsilon,\epsilon]\subset\rho(E)$, each $\xi\in
[-\epsilon,\epsilon]$ is a regular value for~$\rho$ and
for~$\rho\restrict{M\sm C}$, $\rho<-\epsilon$ on $\partial E$,
$\rho>\epsilon$ on~$X\sm W$, and $|\rho|>\epsilon$ on $C\sm Z$. To
see this, we fix a \rel \cpt \nbdns~$V$ of~$\partial E$ in~$X$, we
let $A$ be the union of all the \rel \cpt (in $X$) \concomps of~$C$ whose
closures either meet~$\partial E$ or are contained in~$V$, and we
let $B=C\sm(A\cup Z)$. The closed set $A\cup\partial E$ is \cptns.
For only finitely many \concomps of~$C$ can meet~$\partial V$ and,
of the other \concompsns, those contained in~$A$ must then be
contained in~$V$. The set~$B$ must be closed in $X$, since each of
the finitely many \concomps of~$B$ which meet~$\partial V$ is
closed in~$X$ ($A\cup Z$ contains all of the non-closed \concomps
of~$C$) and the remaining \concomps are contained in
$X\sm\overline V$ ($A$ contains all of the \concomps contained
in~$V$). Thus we may choose a \rel \cpt \nbdns~$W$ of
$A\cup\partial E$ in $X\sm B$ and a nonnegative \cinf \fn $\alpha$
on~$X$ \st $\alpha\equiv 0$ on $A\cup\partial E$, $\alpha>1$ on
$X\sm W$, and $[0,1]\subset\alpha(E)$. In particular, we have
$\alpha\inv((0,1))\cap C\cap M\subset W\cap C\cap M\subset Z\cap M=\emptyset$.
Choosing a number~$s\in
(0,1)$ which is a regular value of both~$\alpha$ and
$\alpha\restrict{(M\sm C)}$, we see that the \nbdns~$W$, the
\fnns~$\rho=\alpha-s$, and any sufficiently small~$\epsilon>0$
have the required properties.

Let $H$ be the \cpt subset of the manifold $M\sm C$ given by
\[
H=h\inv([0,r])\cap\rho\inv([-\epsilon,\epsilon]),
\]
and, for each $\xi\in[-\epsilon,\epsilon]$, let $H_\xi$ be the
fiber $H\cap\rho\inv(\xi)$; a \cpt subset of the submanifold
$\rho\inv(\xi)\cap M=\rho\inv(\xi)\cap(M\sm C)$ of $M\sm C$. There is then a constant
$u>0$ giving the uniform bound
\[
\left|\int_{H_\xi}\eta^{n-1}\restrict{\rho\inv(\xi)\cap (M\sm
C)}\right|<u\qquad\forall\,\xi\in[-\epsilon,\epsilon]
\]
(to see this, one first obtains the local inequality by
considering local coordinates of the form
$(\rho,x_2,\dots,x_{2n-1})$ in $M\sm C$ and then covers $H$ by
finitely many such coordinate \nbdsns). For each $\xi\in[-\epsilon,\epsilon]$,
let $\Omega_\xi=\setof{x\in X}{\rho(x)<\xi}$; a \rel \cpt \nbd of $\partial E$ in~$X$.

Given $a$ and $b$ with $0<a<b<r$, we may choose a regular value
$\xi=\xi(a,b)\in (-\epsilon,\epsilon)$ for
$\rho\restrict{h\inv(a)\cup h\inv(b)}$. The set
\(D=h\inv\left((a,b)\right)\sm\overline\Omega_\xi\) is then a \rel
\cpt open subset of~$M\sm C$ with boundary
\[
\partial D=\overline\Gamma_0\cup\bar\Lambda\cup\overline\Gamma_1;
\]
where $\Gamma_0=h\inv(a)\sm\overline\Omega_\xi$,
$\Gamma_1=h\inv(b)\sm\overline\Omega_\xi$, and
$\Lambda=h\inv\left((a,b)\right)\cap\partial\Omega_\xi$. For if
$t\in(0,r)$, then $(\partial E)\cup h\inv(t)\subset U\Subset X$
for some open set $U$ and $h\inv(I)\cap\partial U=\emptyset$ for
some open interval~$I$ with $t\in I\subset(0,r)$. On the other
hand, each nonempty level of $h$ over a point in~$I$ is both \rel
\cpt in~$X$ and non\cptns, and, therefore, must meet~$U$. Hence
$h\inv(I)\subset U$ and it follows that $D\Subset M$. We have
$\overline D\subset h\inv([a,b])\subset M\sm C$, so $D\Subset M\sm
C$. Moreover $D$ is smooth at each boundary point not in the
``corners" $\overline\Gamma_0\cap\bar\Lambda$ and
$\bar\Lambda\cap\overline\Gamma_1$. Applying Stokes' theorem to
the closed form $\eta^{n-1}$, we get
\[
\int_{\Gamma_1}\eta^{n-1}-\int_{\Gamma_0}\eta^{n-1}=\int_{\Lambda}\eta^{n-1}
\]
(where we take the orientations associated to the \cpx structure
on $\Gamma _0$ and $\Gamma _1$ and the orientation outward from
$\Omega_\xi$ on~$\Lambda$, and we have used the choice of~$\xi$ as a
regular value for $\rho\restrict{h\inv(a)\cup h\inv(b)}$).
On the other hand, we have $\Lambda\subset H_\xi$, so the absolute
value of the integral on the right-hand side of the above equality
is bounded above by the constant~$u$ which does not depend on the
choice of $a$ and $b$. Thus
\[
0\leq\int_{h\inv(a)\sm\overline\Omega_\epsilon}\eta^{n-1}
\leq\int_{\Gamma_0}\eta^{n-1}\leq u+\int_{\Gamma_1}\eta^{n-1}
\leq u+\int_{h\inv(b)\sm\overline\Omega_{-\epsilon}}\eta^{n-1}
\]
and, similarly,
\[
0\leq\int_{h\inv(b)\sm\overline\Omega_\epsilon}\eta^{n-1}
\leq u+\int_{h\inv(a)\sm\overline\Omega_{-\epsilon}}\eta^{n-1}.
\]
Fixing $c\in (0,r)$ and setting $\Omega=\Omega_\epsilon$, we see that
\[
\int_{h\inv(t)\sm\overline\Omega}\eta^{n-1}\leq R\qquad\forall\, t\in(0,r),
\]
where
\[
R=u+\int_{h\inv(c)\sm\overline\Omega_{-\epsilon}}\eta^{n-1}.
\]
The lemma now follows.
\end{pf}

\begin{lem}\label{Rel cpt analytic set in bdd geom end lem}
Let $(X,g)$ be a \con non\cpt complete Hermitian \mfldns, let $E$
be a special end of type~(BG) in~$X$, and let
$A$ be an \anal subset of~$E$ with positive dimension at each point.
Then $A\Subset X$ if and only if
$A\sm\overline\Omega$ has finite volume for some (and, therefore,
for every) \rel \cpt \nbdns~$\Omega$ of~$\partial E$ in~$X$.
\end{lem}
\begin{pf} Suppose $A$ is an \anal subset of~$E$
with $\dim_xA>0$ for each $x\in A$. We proceed as in the first part of the proof of
Theorem~\ref{L2 Castelnuovo de Franchis thm from intro}. Given a
\rel \cpt \nbdns~$\Omega$ of $\partial E$ in $X$, Lelong's
monotonicity formula (see 15.3 in \cite{Chirka}) implies that
there is a constant $c>0$ \st each point $p\in E\sm\Omega$ has a
\nbd $U_p\Subset E$ for which $\diam(U_p)<1$ and $\Vol(D\cap
U_p)\geq c$ for every \con \anal set $D$ of positive dimension in $E$
with $p\in D$. Therefore, if $A\sm\overline\Omega$ has finite
volume, then $\overline A$ must be \cptns. The claim now follows
easily.
\end{pf}

\begin{lem}\label{Rel cpt levels in bdd geom end lem}
Let $(X,g)$ be a \con non\cpt complete Hermitian \mfldns, let $E$
be a special end of type~(BG) in~$X$, and let
$h\colon E\to\mathbb P^1$ be a nonconstant \holo mapping \st the
levels over almost every point in~$\mathbb P^1$ have finite volume
and every nonempty level is non\cptns.
\begin{enumerate}

\item[(a)] If $F$ is a fiber of~$h$ with non\cpt closure in~$X$,
then $F$ has a \concomp~$L$ \st $\overline L$ is non\cpt and
$\overline L\cap\partial E\neq\emptyset$.

\item[(b)] The set \(Q=\setof{\zeta\in\mathbb P^1}{h\inv(\zeta
)\Subset X}=\setof{\zeta\in\mathbb P^1}{L\Subset X\text{ for each
level } L\text{ over }\zeta}\) (the second equality follows
from~(a)) is open and the inverse image of any \cpt subset of~$Q$
is \rel \cpt in~$X$.

\end{enumerate}
\end{lem}
\begin{pf}
For the proof of~(a), suppose $\zeta _0\in\mathbb P^1$ is a point
for which the corresponding fiber $F=h\inv(\zeta_0)$ is not \rel
\cpt in~$X$. Fix a \rel \cpt \nbdns~$\Omega$ of~$\partial E$
in~$X$ and let $L_1,\dots,L_m$ be the (finitely many) \concomps
of~$F$ which meet $\partial\Omega$. Since $h$ has no nonempty \cpt
levels, the closure of any \concomp of~$F$ which does not meet
$\partial\Omega$ must either lie in $\Omega$ and meet~$\partial E$
or lie in $E\sm\overline\Omega$. Thus the union~$H$ of all of the
\concomps of~$F$ which are closed in~$X$ is itself a closed set
in~$X$ which is contained in $L_1\cup\cdots\cup L_m\cup
(E\sm\overline\Omega)$. Therefore, by replacing $\Omega$ with a
\rel \cpt \nbd of $\partial E$ in $X\sm H$, we may assume that
$\overline L_i\cap
\partial E\neq\emptyset$ for $i=1,\dots,m$.

Now if the levels $L_1,\dots,L_m$ are \rel \cpt in~$X$, then, by
replacing $\Omega$ with a \rel \cpt \nbd of the \cpt set
$\overline\Omega\cup L_1\cup\dots\cup L_m$ in $X\setminus H$, we
may assume that $F\cap\partial\Omega=\emptyset$; that is,
$\zeta_0$ lies in the complement $V\equiv\mathbb P^1\sm
h(E\cap\partial\Omega)$ of the \cpt set $h(E\cap\partial\Omega)$.
On the other hand, $F$ meets $E\sm\overline\Omega$ and $h$ is an
open mapping, so $h(E\sm\overline\Omega)$ is a \nbd of $\zeta _0$
in $\mathbb P^1$. Since the levels of~$h$ over almost every point
in~$\mathbb P^1$ have finite volume, Lemma~\ref{Rel cpt analytic
set in bdd geom end lem} implies that there exists a point
$\zeta\in V$ \st $h\inv(\zeta)$ has a \rel \cpt \concompns~$L$
which meets $E\sm\overline\Omega$. Since $h$ has no nonempty \cpt levels,
$L$ must meet $\partial\Omega$ and we have arrived at a
contradiction. Thus there is a \concomp of $F$ which is not \rel \cpt and which is not
closed in~$X$, and~(a) is proved.

For the proof of~(b), given a point \(\zeta _0\in
Q=\setof{\zeta\in\mathbb P^1}{h\inv(\zeta)\Subset X}\), we may
choose a \rel \cpt \nbdns~$\Omega$ of the \cpt set
$h\inv(\zeta_0)\cup\partial E$ in~$X$. The open set $V=\mathbb
P^1\sm h(E\cap\partial\Omega)$ is then a \nbd of $\zeta_0$. If
$\zeta\in\mathbb P^1$ and $F=h\inv(\zeta)$ meets $E\sm\Omega$,
then~(a) implies that $F$ has a (non\cptns) \concompns~$L$ \st
either $\overline L$ is non\cpt and meets~$\partial E$ or $L$ is
\rel \cpt in~$X$ and meets $E\sm\Omega$. In either case, $L$ meets
both $\Omega $ and $E\sm\Omega$, and hence $L$ meets
$\partial\Omega$. It follows that $h\inv(V)\subset\Omega\Subset
X$. In particular, $V\subset Q$ and hence $Q$ is open. Moreover if
$K\subset Q$ is a \cpt subset, then one may cover $K$ by finitely
many such \nbdsns~$V$ with \rel \cpt inverse image and, therefore,
$h\inv(K)\Subset X$. Thus~(b) is proved.
\end{pf}
\begin{lem}\label{Rel cpt levels over regular values in bdd geom end lem}
Suppose $(X,g)$ is a \con non\cpt complete Hermitian \mfldns, $E$
is a special end of type~(BG) in~$X$, $h\colon E\to\mathbb P^1$ is a nonconstant \holo
mapping, $C$ is the set of critical points of~$h$, and $Z$ is the
union of all \concompsns~$A$ of~$C$ for which $\bar A$ is non\cpt
and $\bar A\cap\partial E\neq\emptyset$. Assume that $g\restrict E$ is K\"ahler, the levels
of~$h$ over almost every point in~$\mathbb P^1$ have finite
volume, and every nonempty level is non\cptns. Then the fiber of
$h$ over every point in~$\mathbb P^1\sm h(Z)$ is \rel \cpt in~$X$.
\end{lem}
\begin{pf}
We must show that each point $\zeta_0\in\mathbb P^1\sm h(Z)$ lies
in the open set
\[
Q\equiv\setof{\zeta\in\mathbb P^1}{h\inv(\zeta )\Subset X}.
\]
As in the proof of Theorem~\ref{L2 Castelnuovo de Franchis thm
from intro}, the idea is to consider converging sequences of
fibers with uniformly bounded volume (as in \cite{Gro-Kahler
hyperbolicity}, \cite{ArapuraBressRam}, and \cite{ABCKT}). We
first observe that, by Lemma~\ref{Rel cpt analytic set in bdd geom
end lem} and Lemma~\ref{Rel cpt levels in bdd geom end lem},
$\mathbb P^1\sm Q$ is a closed set of measure~$0$. Therefore,
since the set of critical values of~$h$ is countable, we may
choose an extended real line~$\ell$ through $\zeta_0$ \st
$\ell\sm\set{\zeta_0}\subset\mathbb P^1\sm h(C)$ and $\ell\sm Q$
is a (closed) set of (real $1$-dimensional) measure~$0$. If $\ell$
is not contained in~$Q$, then there is an open segment $I$
in~$\ell\cap Q$ which has a boundary point~$a\in\ell$ and another boundary point
not equal to~$a$. Fixing a sufficiently large \rel \cpt \nbdns~$\Omega$
of~$\partial E$ in~$X$ and applying Lemma~\ref{finite volume
fibers in end lem}, we get a constant~$v$ \st
\[
\Vol(h\inv(t)\sm\overline\Omega)<v\quad\forall\,t\in I.
\]
As in the proof of~Lemma~\ref{Rel cpt analytic set in bdd geom end
lem}, this volume bound also gives a bound on the diameter
relative to the distance \fn on~$X$. For Lelong's monotonicity
formula (see 15.3 in \cite{Chirka}) implies that there is a
constant $c>0$ \st each point $p\in E\sm\Omega$ has a \nbd
$U_p\Subset E$ for which $\diam(U_p)<1/2$ and $\Vol(D\cap U_p)\geq
c$ for every \con \anal set $D$ of positive dimension in $E$ with $p\in
D$. We may also fix a point $O\in X$ and a constant~$R>0$ \st
$\Omega\subset B(O;R)$. If some level $L$ of $h$ over a point
$t\in I\subset Q$ meets $E\sm B(O;R+k)$ for some positive
integer~$k$, then, since $L$ also meets $\Omega$, $L$ meets
$\partial B(O;R+j)$ for all $j=1,\dots,k$. The volume estimate
then implies that
\[
v>\Vol(L\sm\overline\Omega)\geq kc.
\]
Thus $h\inv(t)\subset B=B(O;R+(v/c)+1)\Subset X$ for all $t\in I$.
But $h\inv(a)$ is not \rel \cpt in~$X$ and $h$ is an open mapping,
so $h(E\sm\overline B)$ is an open set which contains~$a$ and,
therefore, a point $t\in I$. Thus we have arrived at a
contradiction and hence $\zeta_0\in\ell\subset Q$.
\end{pf}
\begin{pf*}{Proof of Theorem~\ref{L2 Castelnuovo for end thm}}
It suffices to show that the nonconstant \holo mapping
\[
f\equiv\frac{\omega_1}{\omega_2}\colon E\to\mathbb P^1
\]
has a nonempty \cpt level. Assuming $f$ has no nonempty \cpt levels, we will apply
the above lemmas to obtain a contradiction. For this, we let~$C$
be the set of critical points and we let $Z$ be the union of all
\concomps~$A$ of~$C$ for which $\bar A$ is non\cpt and
meets~$\partial E$.

According to Lemma~\ref{Finite volume for level lemma}, the levels
over almost every point in~$\mathbb P^1$ have finite volume.
According to Lemma~\ref{Rel cpt levels in bdd geom end lem} and
Lemma~\ref{Rel cpt levels over regular values in bdd geom end
lem}, the set
\[
Q=\setof{\zeta\in\mathbb P^1}{f\inv(\zeta )\Subset
X}=\setof{\zeta\in\mathbb P^1}{L\Subset X\text{ for each level }
L\text{ over }\zeta}
\]
is an open set containing $\mathbb P^1\sm f(Z)$. Moreover, $Z$ has
only finitely many \concompsns, since each \concomp must meet the
(\cptns) boundary of any fixed \rel \cpt \nbd of~$\partial
E$ in~$X$. Consequently, $f(Z)\supset\mathbb P^1\sm Q$ is a finite
set and Lemma~\ref{Rel cpt levels in bdd geom end lem} implies
that $Q\neq\mathbb P^1$. Thus $\mathbb P^1\sm
Q=\set{\zeta_1,\dots,\zeta_m}$ for distinct points
$\zeta_1,\dots,\zeta_m$ and we may choose domains $U_1,\dots,U_m$
in $\mathbb P^1$ \stns, for each $i=1,\dots,m$, we have
$\zeta_i\in U_i$, $\overline{U_i}\cap\overline{U_j}=\emptyset$
for all $j\neq i$, $0\in U_i$ if and only if $\zeta_i=0$,
and $\infty\in U_i$ if and only if $\zeta_i=\infty$. We may fix a \rel \cpt \nbdns~$\Omega$
of~$\partial E$ in~$X$ so that $\Omega$ contains $f\inv(\infty)$
if $\infty\in Q$.  Lemma~\ref{Rel cpt levels in bdd geom end lem}
implies that the set $H=\overline\Omega\cup f\inv\left(\mathbb
P^1\sm(U_1\cup\cdots\cup U_m)\right)$ is \cptns. Let $E_0$ be a
\concomp of~$X\sm H$ which is contained in $E$ and which has
non\cpt closure in~$X$. We have $f(E_0)\subset U_i$ for some~$i$
and we may define the \holo \fn $h\colon E_0\to\C$ by
\[
h=\left\{
\begin{alignedat}{2}
&1/f\restrict{E_0}
&\quad\text{ if }\zeta_i=\infty\\
&f\restrict{E_0}-\zeta_i&\quad\text{ if }\zeta_i\in \C
\end{alignedat}
\right.
\]
The inverse image of any \cpt subset of $\overline
U_i\sm\set{\zeta_i}$ is \rel \cpt in $X$, so $h$ must vanish at
infinity in~$X$. Proposition~\ref{parabolic end vanish at infinity
end prop} then implies that some nonempty open subset of~$E_0$
admits a proper \holo mapping onto a Riemann surface. In
particular, $f$ has a \cpt level in~$E_0\subset E$ and we have
arrived at a contradiction. Thus the theorem is proved.
\end{pf*}
\begin{rmk}
If we assume that the end $E$ is hyperbolic in~$(X,g)$, then we
need only apply Lemma~\ref{holo fn vanish at infinity parabolic
lemma} in place of Proposition~\ref{parabolic end vanish at
infinity end prop} in order to obtain a contradiction. So the
proof in this case (which is the case required for~\cite{NR-BH
bound geom hyperbolic Kahler}) is simpler.
\end{rmk}

\bibliographystyle{amsalpha.bst}

\end{document}